\documentclass[fleqn]{article}
\usepackage{amssymb,latexsym,theorem}

\newcommand{\G}{\mathcal{G}}

\newcommand{\B}{\mathbf{B}}

\newcommand{\bC}{\mathbb{C}}

\newcommand{\bZ}{\mathbb{Z}}

\newcommand{\PG}{\mathrm{PG}}

\newcommand{\F}{\mathbb{F}}

\newtheorem{lemma}{Lemma}[section]
\newtheorem{prop}[lemma]{Proposition}
\newtheorem{theo}[lemma]{Theorem}
\newtheorem{co}[lemma]{Corollary}
\newtheorem{conj}[lemma]{Conjecture}

\def\pr{\noindent{\bf Proof. }}
\def\eop{\hspace*{\fill}$\Box$}
\title{On certain submodules of Weyl modules for $\mathrm{SO}(2n+1,\mathbb{F})$ with $\mathrm{char}(\mathbb{F}) = 2$}
\author{Ilaria Cardinali and Antonio Pasini}
\date{}

\begin{document}
\maketitle
\begin{abstract}
For $k = 1, 2,...,n-1$ let $V_k = V(\lambda_k)$ be the Weyl module for the special orthogonal group $G = \mathrm{SO}(2n+1,\F)$ with respect to the $k$-th fundamental dominant weight $\lambda_k$ of the root system of type $B_n$ and put $V_n = V(2\lambda_n)$. It is well known that all of these modules are irreducible when $\mathrm{char}(\F) \neq 2$ while when $\mathrm{char}(\F) = 2$ they admit many proper submodules. In this paper, assuming that $\mathrm{char}(\F) = 2$, we prove that $V_k$ admits a chain of submodules $V_k = M_k \supset M_{k-1}\supset ... \supset M_1\supset M_0 \supset M_{-1} = 0$ where $M_i \cong V_i$ for $1,..., k-1$ and $M_0$ is the trivial 1-dimensional module. We also show that for $i = 1, 2,..., k$ the quotient $M_i/M_{i-2}$ is isomorphic to the so called $i$-th Grassmann module for $G$. Resting on this fact we can give a geometric description of $M_{i-1}/M_{i-2}$ as a submodule of the $i$-th Grassmann module. When $\F$ is perfect $G\cong \mathrm{Sp}(2n,\F)$ and $M_i/M_{i-1}$ is isomorphic to the Weyl module for $\mathrm{Sp}(2n,\F)$ relative to the $i$-th fundamental dominant weight of the root system of type $C_n$.  All irreducible sections of the latter modules are known. Thus, when $\F$ is perfect, all irreducible sections of $V_k$ are known as well.
 \end{abstract}

{\scriptsize
}

\section{Introduction}\label{Introduction}

Let $V := V(2n+1,\F)$ be a $(2n+1)$-dimensional vector space over a field $\F$ and let $\eta$ be a non singular quadratic form of $V$ of Witt index $n$. Let $\Delta$ be the building of type $B_n$ where the elements of type $k = 1, 2,..., n$ are the $k$-dimensional subspaces of $V$ totally singular for $\eta$.

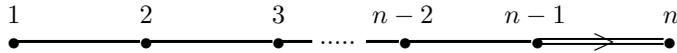
\begin{picture}(310,36)(0,0)
\put(20,8){$\bullet$}
\put(23,11){\line(1,0){47}}
\put(70,8){$\bullet$}
\put(73,11){\line(1,0){47}}
\put(120,8){$\bullet$}
\put(123,11){\line(1,0){12}}
\put(138,10){$.....$}
\put(156,11){\line(1,0){12}}
\put(168,8){$\bullet$}
\put(171,11){\line(1,0){47}}
\put(218,8){$\bullet$}
\put(241,8){\large{$>$}}
\put(221,10){\line(1,0){47}}
\put(221,12){\line(1,0){47}}
\put(268,8){$\bullet$}
\put(20,18){1}
\put(70,18){2}
\put(120,18){3}
\put(158,18){$n-2$}
\put(208,18){$n-1$}
\put(268,18){$n$}
\end{picture}

\noindent
For $1\leq k\leq n$, let $\G_{k}$ be the $k$-grassmannian of $\PG(V)$, that is $\G_{k}$ is the point-line geometry where the points are the $k$-dimensional subspaces of $V$ and the lines are the sets $l_{X,Y} = \{Z\mid X \subset Z\subset Y , ~\mathrm{dim}(Z) = k\}$ for a given pair $(X,Y)$ of subspaces of $V$ where $\mathrm{dim}(X) = k-1$ and  $\mathrm{dim}(Y) = k+1$. (Note that $X = 0$ when $k = 1$.)

Let $\Delta_{k}$ be the $k$-{\em grassmannian}  of $\Delta$ The point-line geometry $\Delta_{k}$ is the proper subgeometry of $\G_{k}$ defined as follows. The points of $\Delta_{k}$ are the $k$-elements of $\Delta$. When $1\leq k < n$ the lines of $\Delta_{k}$ are the lines $l_{X,Y}$ of $\G_{k}$ with  $Y$ totally singular. Note that $\Delta_{1}$ is the polar space associated to $\eta$.

Let $k=n$. Then the lines of $\Delta_{n}$ are the sets
\begin{equation}\label{lines}
l_{X} = \{Z\mid X \subset Z\subset X^{\perp} , ~\mathrm{dim}(Z) = n, ~ Z ~\mbox{totally singular}\}
\end{equation}
where $X$ is a given $(n-1)$-dimensional totally singular subspace of $V$ and $X^\perp$ is the orthogonal of $X$ with respect to (the linearization of) $\eta$. The vector space $X^\perp/X$ is 3-dimensional and $l_{X}$ is a non-singular conic in the projective plane $\PG(X^\perp/X)$. The geometry $\Delta_{n}$ is called the {\em dual polar space} of type $B_n$.

Let $1\leq k\leq n$ and $W_{k}:=\bigwedge^k V$. Then $\mathrm{dim}(W_{k}) = {{2n+1}\choose k}$.
Let $\iota_{k}\colon \G_{k}\rightarrow W_{k}$ be the embedding of $\G_{k}$ in $\mathrm{PG}(W_{k})$ mapping every  point of $\G_{k}$, i.e. every $k$-dimensional subspace $\langle v_1,v_2,\dots, v_k\rangle$ of $V$, onto the point $\langle v_1\wedge v_2\wedge\dots \wedge v_k \rangle$ of $\PG(W_{k})$. Let $\varepsilon_{k}:= \iota_{k}|_{\Delta_{k}}$ be the restriction of $\iota_{k}$ to $\Delta_{k}$. The map $\varepsilon_{k}$ is called the {\em grassmann embedding} of $\Delta_{k}$. If $k<n$ then ${\varepsilon}_{k}$ maps lines of $\Delta_{k}$ onto projective lines. If $k=n$ then ${\varepsilon}_{k}$ maps lines of $\Delta_{k}$ onto non-singular conics of $W_{k}$.

Put $W_{k}^{\circ}:=\langle \varepsilon_{k}(\Delta_{k})\rangle$. As proved in Cardinali and Pasini \cite{CP1},
\begin{equation}\label{dim-grass}
\mathrm{dim}(W^\circ_{k}) = \left\{\begin{array}{ll}
{{2n+1}\choose k} & \mbox{if}~ \mathrm{char}(\F) \neq 2,\\
 & \\
{{2n+1}\choose k}-{{2n+1}\choose{k-2}} & \mbox{if}~\mathrm{char}(\F) = 2.
\end{array}\right.
\end{equation}
So, if either $\mathrm{char}(\F)\neq 2$ or $k = 1$ then $W^{\circ}_{k} = W_{k}$  (recall that ${{2n+1}\choose{-1}} := 0$) while if $\mathrm{char}(\F) = 2$ and $k > 1$ then $W^{\circ}_{k}$ is a proper subspace of $W_{k}$ of codimension ${{2n+1}\choose{k-2}}$.

Let $G := \mathrm{SO}(2n+1,\F)$ be the stabilizer of the form $\eta$ in the special linear group $\mathrm{SL}(V) = \mathrm{SL}(2n+1,\F)$. The group $G$ also acts on $W_{k}$, according to the following rule:
\[g(v_1\wedge...\wedge v_k) = g(v_1)\wedge...\wedge g(v_k) ~~ \mbox{for $g\in G$ and $v_1,..., v_k\in V$}.\]
Note that $\mathrm{SO}(2n+1,\F) = \mathrm{PSO}(2n+1,\F)$, namely $G$ is the adjoint Chevalley group of type $B_n$ defined over $\F$. The universal Chevalley group of type $B_n$ is the spin group $\widetilde{G} = \mathrm{Spin}(2n+1,\F)$.
If $\mathrm{char}(\F) = 2$ then $\widetilde{G} = G$. If $\mathrm{char}(\F) \neq 2$ then $\widetilde{G} = 2^\cdot G$, a non-split central extension of $G$ by a group of order two.

Clearly, $W_{k}^{\circ}$ is a $G$-module. It is called the {\em $k$-grassmann orthogonal module} for $G$  (also {\em grassmann orthogonal module} if $k$ is clear from the context).

Let $\lambda_1, \lambda_2,..., \lambda_n$ be the fundamental dominant weights for the root system of type $B_n$, the nodes of the $B_n$-diagram being numbered in the usual way (see the picture at the beginning of this introduction). For a positive integral combination $\lambda$ of $\lambda_1, \lambda_2,..., \lambda_n$ we denote by $V(\lambda)$ the Weyl module for $\widetilde{G}$ with $\lambda$ as the highest weight. If $\widetilde{G}$ acts unfaithfully on $V(\lambda)$ then it induces $G$ on $V(\lambda)$, namely $V(\lambda)$ is also a $G$-module. For instance, $V(\lambda_k)$ is a $G$-module for every $k < n$ as well as $V(2\lambda_n)$ while if $\mathrm{char}(\F)\neq 2$ then $V(\lambda_n)$ is a $\widetilde{G}$-module but not a $G$-module.

We introduce the following notation not to be forced to split our statements in two distinct cases, namely $k < n$ and $k = n$. We put $\lambda^\circ_k = \lambda_k$ for $k < n$ and $\lambda^\circ_n = 2\lambda_n$. We also put $\lambda^\circ_0 = 0$ (the null weight). Thus $V(\lambda_0^\circ)$ is the trivial ($1$-dimensional) $G$-module.
Note that
\begin{equation}\label{dim-Weyl}
\mathrm{dim}(V(\lambda_k^\circ)) = {{2n+1}\choose k}
\end{equation}
for every $k = 0, 1,..., n$, as one can see by applying the Weyl dimension formula (see \cite[24.3]{Humphreys}, for instance). The following is well known (see e.g. Blok~\cite{Blok2011}; see also Theorem \ref{immagini} of Section \ref{general setting} of this paper).

\begin{theo}\label{theorem grassmann-weyl}
For $k = 1, 2,..., n$ the grassmann orthogonal module $W_{k}^{\circ}$ is a homomorphic image of the Weyl module $V(\lambda_k^\circ)$.
\end{theo}
Therefore, in view of (\ref{dim-grass}) and (\ref{dim-Weyl}),

\begin{co}\label{corollary grassmann-weyl}
If either $\mathrm{char}(\F) \neq 2$ or $k = 1$ then $W^\circ_{k} \cong V(\lambda^\circ_k)$. If $\mathrm{char}(\F) = 2$ and $k > 1$ then $W^\circ_{k} $ is a proper homomorphic image of $V(\lambda^\circ_k)$.
\end{co}

It is worth spending a few more words on the statement of Theorem \ref{theorem grassmann-weyl}. We recall that $V(\lambda^\circ_k)$ is primarily defined as a module for the enveloping algebra of the Lie algebra of $G$. A choice of a $BN$-pair of $G$, equivalently an apartment $\Sigma$ of $\Delta$ and a chamber $C$ of $\Sigma$, is implicit in this way of defining $V(\lambda^\circ_k)$ in such a way that if $V^+$ is the (1-dimensional) weight space of $V(\lambda^\circ_k)$ of weight $\lambda^\circ_k$ and $X_k := \langle e_1,..., e_k\rangle$ is the $k$-element of $C$ then the parabolic subgroup of $G$ stabilizing $X_k$ is the stabilizer of $V^+$ in $G$. The homomorphism $\varphi_k$ from $V(\lambda^\circ_k)$ to $W^\circ_{k}$ maps $V^+$ onto $\langle e_1\wedge e_2\wedge...\wedge e_k\rangle$. This condition uniquely determines $\varphi_k$ modulo scalars. We denote the kernel of $\varphi_k$ by $K_k$. Thus, $W^\circ_{k}\cong V(\lambda^\circ_k)/K_k$. In view of Corollary \ref{corollary grassmann-weyl} we have $K_k = 0$ when either $\mathrm{char}(\F) \neq 2$ or $k = 1$.

Let $\mathrm{char}(\F) = 2$. Then the bilinear form associated to $\eta$ is degenerate with $1$-dimensional radical $N_0$. The subspace $N_0$ is called the {\em nucleus} of the quadratic form $\eta$. It is the unique $1$-dimensional subspace of $V$ stabilized by $G$ and it is non-singular.

Suppose firstly that  $k > 1$. Given an element $X$ of $\Delta$ of type $k-1$, let $\mathrm{St}(X)$ be its upper residue, formed by the elements of $\Delta$ of type $k, k+1,..., n$ that contain $X$. We call $\mathrm{St}(X)$ the {\em star} of $X$. Clearly, $\mathrm{St}(X)$ is (the building of) an orthogonal polar space of rank $n-k+1$ defined in $X^\perp/X$. Let $n_X$ be the nucleus of a quadratic form of $X^\perp/X$ associated to the polar space $\mathrm{St}(X)$. Then $n_X = N_X/X$ where $N_X = \langle X, N_0\rangle$. Hence $N_X$ is a point of $\G_{k}$ and, since $n_X$ belongs to $X^\perp/X$, which is spanned by the 1-dimensional subspaces $Y/X$ for $Y$ ranging in the set of points of $\mathrm{St}(X)$, the point $\iota_{k}(N_X)$ of $\PG(W_{k})$ belongs to $\mathrm{PG}(W^\circ_{k})$. We define the $k$-{\em nucleus subspace} $\overline{\cal N}_k$ of $W^\circ_{k}$ as the subspace of $W^\circ_{k}$ spanned by the 1-dimensional subspaces $\iota_{k}(N_X)$ for $X\in \Delta_{k-1}$ and we denote its preimage in $V(\lambda^\circ_k)$ by the symbol ${\cal N}_k$. Thus, ${\cal N}_k/K_k\cong\overline{\cal N}_k$. Note that, since the homomorphism $\varphi_k:V(\lambda^\circ_k)\rightarrow W^\circ_{k}$ is uniquely determined modulo scalars and ${\cal N}_k = \varphi^{-1}_k(\overline{\cal N}_k)$, the submodule ${\cal N}_k$ of $V(\lambda_k^\circ)$ is uniquely determined as well.

When $k = 1$ we put ${\cal N}_1 = N_0$, regarded as a subspace of $V(\lambda^\circ_1)$ (as we may, since $K_1 = 0$). Clearly, $\mathcal{N}_{k}$ is stabilized by $G$, whence it is a submodule of $V(\lambda^\circ_k)$. The following is proved in Cardinali and Pasini~\cite{CP1}.

\begin{theo}\label{theorem 1}
Let $\mathrm{char}(\F) = 2$. Then $\mathrm{dim}({\cal N}_k) = \mathrm{dim}(V(\lambda_{k-1}^\circ)) = {{2n+1}\choose{k-1}}$ for every $k = 1, 2,..., n$, namely ${\cal N}_k$ has codimension ${{2n+1}\choose k}-{{2n+1}\choose{k-1}} = {{2n}\choose k}-{{2n}\choose{k-2}}$ in $V(\lambda^\circ_k)$.
\end{theo}

The following theorem, to be proved in Section \ref{section-main theorem}, is the main result of this paper.

\begin{theo}\label{main theorem}
Let $\mathrm{char}(\F) = 2$. Then ${\cal N}_k \cong V(\lambda_{k-1}^\circ)$ (isomorphism of $G$-modules) for every $k = 1, 2,..., n$. Moreover, if $k > 1$ then the isomorphism from ${\cal N}_k$ to $V(\lambda_{k-1}^\circ)$ can be chosen in such a way that it maps $K_k$ onto the submodule ${\cal N}_{k-1}$ of $V(\lambda_{k-1}^\circ)$.
\end{theo}

The next corollaries immediately follow from the previous theorem.

\begin{co}\label{main corollary}
Let $\mathrm{char}(\F) = 2$. Then $K_k \cong V(\lambda_{k-2}^\circ)$ for every $k = 2, 3,..., n$.
\end{co}

\begin{co}\label{main corollary bis}
Let $\mathrm{char}(\F) = 2$ and $1\leq k \leq n$. Then $V(\lambda^\circ_k)$ admits a chain of submodules
$0 \subset M_0 \subset M_1 \subset ... \subset M_{k-1} \subset M_k = V(\lambda^\circ_k)$
with $M_i\cong V(\lambda^\circ_i)$ for $i = 0, 1,..., k-1$. Moreover, $M_i/M_{i-2}\cong W^\circ_{i}$ for every $i = 2, 3,..., k$ .
\end{co}

The statements of Theorem \ref{theorem 1} and Corollary \ref{main corollary bis} can be made sharper when $\F$ is a perfect field of caracteristic $2$. We shall discuss this case in a few lines, but firstly we recall a few well known facts regarding symplectic polar grassmannians and their natural embeddings.

Put $\overline{V} := V(2n,\F)$. Let $\overline{\Delta}$ be the building of type $C_n$ associated with the symplectic group $\overline{G} := \mathrm{Sp}(2n,\F)$ in its natural action on $\overline{V}$ and, for $k = 1, 2,..., n$, let $\overline{\Delta}_{k}$ be the $k$-grassmannian of $\overline{\Delta}$. Then $\overline{\Delta}_{k}$ is a subgeometry of the $k$-grassmannian ${\overline{\G}}_{k}$ of $\PG(\overline{V})$. Put ${\overline{W}}_{k}:=\bigwedge^k \overline{V}$ and let ${\bar{\iota}}_{k}:{\overline{\G}}_{k}\rightarrow \PG({\overline{W}}_{k})$ be the embedding of ${\overline{\G}}_{k}$ mapping every $k$-dimensional subspace $\langle v_1,..., v_k\rangle$ of $\overline{V}$ onto the point $\langle v_1\wedge ... \wedge v_k\rangle$ of $\PG({\overline{W}}_{k})$. Let $\bar{varepsilon}_{k} := \bar{\iota}_{k}|_{\overline{\Delta}_{k}}$ be the restriction of $\bar{\iota}_{k}$ to $\overline{\Delta}_{k}$ and let $\overline{W}^{\mathrm{sp}}_{k}$ be the subspace of $\overline{W}_{k}$ spanned by $\bar{\varepsilon}_{k}(\overline{\Delta}_{k})$. Then $\bar{\varepsilon}_{k}$ is a projective embedding of $\overline{\Delta}_{k}$ in $\mathrm{PG}(\overline{W}^{\mathrm{sp}}_{k})$, called the {\em natural} or {\em Grassmann} embedding of $\overline{\Delta}_{k}$. The subspace $\overline{W}_{k}^{\mathrm{sp}}$ of $\overline{W}_{k}$ is stabilized by $\overline{G}$ in its natural action on $\overline{W}_{k}$. Thus, $\overline{W}^{\mathrm{sp}}_{k}$ is a $\overline{G}$-modue. It is well known that $\overline{W}^{\mathrm{sp}}_{k}$ has codimension ${{2n}\choose{k-2}}$ in $\overline{W}_{k}$ (see e.g. \cite{BCo1}). Moreover, $\overline{W}^{\mathrm{sp}}_{k}$ is isomorphic to the Weyl module for $\overline{G}$ associated to the $k$-th fundamental dominant weight of the root system of type $C_n$. For the sake of completeness we denote by $\overline{W}^{\mathrm{sp}}_{0}$ the trivial (1-dimensional)  $\overline{G}$-module.

Suppose now that $\F$ is a perfect field of characteristic $2$. Then $G \cong\overline{G}$ and $\Delta\cong\overline{\Delta}$. Accordingly, $\overline{W}^{\mathrm{sp}}_{k}$ can be regarded as a $G$-module. The following is proved in Cardinali and Pasini~\cite{CP1}:

\begin{theo}\label{theorem 1 bis}
Let $\F$ be a perfect field of characteristic $2$. Then $V(\lambda^\circ_k)/{\cal N}_k\cong\overline{W}^{\mathrm{sp}}_{k}$ and ${\cal N}_k/K_k\cong \overline{W}^{\mathrm{sp}}_{k-1}$, for every $k = 1, 2,..., n$.
\end{theo}

Consequently,

\begin{co}\label{main corollary ter}
Let $\F$ be a perfect field of characteristic $2$ and $1\leq k \leq n$. Then $V(\lambda^\circ_k)$ admits a chain of submodules
$0 \subset M_0 \subset M_1 \subset ... \subset M_{k-1} \subset M_k = V(\lambda^\circ_k)$
with $M_i\cong V(\lambda^\circ_i)$ and $M_i/M_{i-1}\cong \overline{W}^{\mathrm{sp}}_{i}$ for $i = 1, 2,..., k$.
\end{co}

In general the series $0 \subset M_0 \subset M_1 \subset ... \subset M_{k-1} \subset M_k = V(\lambda^\circ_k)$ is not a composition series. By Corollary \ref{main corollary ter}, when $\F$ is perfect we can produce a composition series of $V(\lambda^\circ_k)$ by inserting a composition series of $\overline{W}^{\mathrm{sp}}_{i}$ between $M_{i-1}$ and $M_i$, for every $i = 1, 2,..., k$. In other words, when $\F$ is perfect the irreducible sections of $V(\lambda^\circ_k)$ are the irreducible sections of $\overline{W}^{\mathrm{sp}}_{i}$ for $i = 1, 2,..., k$. The latters are known (Baranov and Suprunenko \cite{BS}; also Premet and Suprunenko \cite{PS}). Therefore:

\begin{co}\label{final corollary}
Let $\F$ be a perfect field of characteristic $2$ and $1\leq k \leq n$. Then all irreducible sections of  $V(\lambda^\circ_k)$ are known.
\end{co}

It is likely that $V(\lambda^\circ_k)$ admits a unique series of submodules $M_0, M_1,..., M_k$ satisfying the properties of Corollary \ref{main corollary bis}. In Section \ref{uniqueness}, assuming that $\F$ is perfect, we shall prove that this is indeed the case for $k \leq 4$.

\section{Preliminaries}\label{section-associative algebras}

In this section we recall a number of known facts on Lie algebras, Chevalley groups and their Weyl modules, either assumed as well known in Section \ref{Introduction} or to be exploited in the proof of Theorem \ref{main theorem}. We firstly recapitulate the construction of Chevalley groups from simple Lie algebras (Subsection \ref{general setting}). We rely on Humphreys~\cite{Humphreys} and Steinberg~\cite{Stein} for our exposition. We shall give some details, for the commodity of the reader, but not so many. The reader is referred to \cite{Humphreys} and \cite{Stein} for all we shall miss. In the second part of this section (Subsections \ref{orthogonal algebra} and \ref{symplectic algebra}) we shall give more detailed informations on orthogonal and symplectic Lie algebras.

\subsection{General setting}\label{general setting}

Let $\mathfrak{L}_{\mathbb{C}}= \mathfrak{H}\oplus(\oplus_{\alpha\in \Phi}\mathfrak{X}_{\alpha})$ be a simple Lie algebra of rank $n$ over the complex field $\bC$, where $\mathfrak{H}$ is the Cartan subalgebra, $\Phi$ is the set of all roots and $\mathfrak{X}_{\alpha}$ is the $1$-dimensional subalgebra of $\mathfrak{L}_{\mathbb{C}}$ corresponding to the root $\alpha$.
Let $\Pi = \{\alpha_1,..., \alpha_{n}\}$ be a set of simple roots of $\mathfrak{L}_{\mathbb{C}}$ and $\Phi^+$ and $\Phi^-$ the set of roots that are positive or respectively negative with respect to $\Pi$. Let $\{H_1,..., H_{n}\}\cup\{X_\alpha, Y_\alpha\}_{\alpha\in\Phi^+}$ be a Chevalley basis of $\mathfrak{L}_{\bC}$, where $H_i\in \mathfrak{H}$, $X_\alpha\in\mathfrak{X}_{\alpha}$ and $Y_\alpha\in\mathfrak{X}_{-\alpha}$. Denote by  $\mathfrak{A}_{\bC}$ the universal enveloping (associative) algebra of $\mathfrak{L}_{\bC}.$  Chosen an ordering on the set $\Phi^+$, put
\[\begin{array}{lcl}
{\cal A}^- := \{A_t^-\}_{t\in T}, & & A_t^- :=\prod_{\alpha\in\Phi^+} {\frac{Y_{\alpha}^{t_\alpha}}{t_\alpha!}}, \\
{} & & \\
{\cal A}^0 := \{A_s^0\}_{s\in S}, & & A_s^0 :=\prod_{i=1}^{n} {{H_i}\choose{s_i}},\\
{} & & \\
{\cal A}^+ := \{A_r^+\}_{r\in T}, & & A_r^+ :=\prod_{\alpha\in\Phi^+} {\frac{X_{\alpha}^{r_\alpha}}{r_\alpha!}}
\end{array}\]
where
\[{{H_i}\choose{s_i}} := \frac{H_i(H_i- I)...(H_i-(s_i-1)I)}{s_i!}\]
($I$ being the identity element of $\mathfrak{A}_{\bC}$), the factors $Y_\alpha^{t_\alpha}/t_\alpha!$ and $X_\alpha^{r_\alpha}/r_\alpha!$ in $A^-_t$ and $A^+_r$ are ordered according to the ordering previously chosen on $\Phi^+$, the letters $T$ and $S$ denote the set of mappings from $\Phi^+$ and $\{1,2,..., n\}$ respectively to $\mathbb N$ and, given $t, r\in T$ and $s\in S$, we write $t_\alpha$, $r_\alpha$ and $s_i$ for $t(\alpha)$, $r(\alpha)$ and $s(i)$.

The set ${\cal A} := {\cal A}^-\cdot{\cal A}^0\cdot{\cal A}^+ = \{A_t^-A^0_sA^+_r\}_{t,r\in T; s\in S}$ is a linear basis of the algebra $\mathfrak{A}_{\bC}$ (Humphreys~\cite[Chapter 26]{Humphreys}). We call the elements of ${\cal A}^-$, ${\cal A}^0$ and ${\cal A}^+$ {\em monomials} of $(-)$-, $(0)$- and $(+)$-{\em type} respectively and the elements of ${\cal A}$ just {\em monomials}. The identity element $I$ of $\mathfrak{A}_{\bC}$ is regarded as a monomial of $(-)$-type as well as $(0)$- and $(+)$-type.

Every positive root $\alpha \in \Phi^+$ is a sum $\alpha = \sum_{i=1}^nk_{\alpha,i}\alpha_i$ where the coefficients $k_{\alpha,i}$ are nonnegative integers. For $A^+_t\in{\cal A}^+$ and $i = 1, 2,..., n$ we set $d_i(A^+_t) = \sum_{\alpha\in\Phi}t_\alpha k_{\alpha,i}$ and we call $d_i(A^+_t)$ the $i$-{\em degree} of $A^+_t$. The nonnegative integer $d(A^+_t) := \sum_{i=1}^nd_i(A^+_t)$ is the {\em total degree} of $A_t^+$. Similarly, we put $d_i(A^-_t) = -\sum_{\alpha\in\Phi}t_\alpha k_{\alpha,i}$ ($i$-{\em degree} of $A^-_t$) and $d(A^-_t) := \sum_{i=1}^nd_i(A^-_t)$ ({\em total degree} of $A_t^-$).  The $i$-{\em degree} (the {\em total degree}) of a monomial $A^-_tA^0_sA^+_r$  is $d_i(A^-_tA^0_sA^+_r) := d_i(A^-_t)+d_i(A^+_r) = \sum_{\alpha\in\Phi^+}k_{\alpha,i}(r_\alpha-t_\alpha)$ (respectively  $d(A^-_tA^0_sA^+_r) := d(A^-_t)+d(A^+_r)$).

Regarding $\mathfrak{A}_{\bC}$ as a ring, we denote by $\mathfrak{A}^-_{\mathbb{Z}}$ and $\mathfrak{A}_{\mathbb{Z}}$ the subrings of $\mathfrak{A}_{\bC}$ generated by the monomials of type $(-)$ and all monomials respectively. The ring $\mathfrak{A}^-_{\mathbb{Z}}$, regarded as a $\mathbb{Z}$-module, is free over the basis ${\cal A}^-$ while $\mathfrak{A}_{\mathbb{Z}}$ is free over $\cal A$. We also put $\mathfrak{L}_\mathbb{Z} := \mathfrak{L}_{\mathbb{C}}\cap \mathfrak{A}_{\mathbb{Z}} = \langle \{H_i\}_{i=1}^n\cup\{X_\alpha\}_{\alpha\in\Phi}\rangle_{\mathbb{Z}}$.

Let $\Lambda$ be the set of weights of the root system $\Phi$ and $\Lambda^+$ the set of dominant weights (relative to the choice of $\Pi = \{\alpha_1,...\alpha_n\}$ as the set of simple roots). As in Humphreys \cite{Humphreys}, for two weights $\lambda, \mu \in \Lambda$ we write $\mu \prec \lambda$ if $\lambda-\mu$ is a sum of positive roots.

Given a dominant weight $\lambda\in \Lambda^+$ let $\bf{v}^+$ be a vector of highest weight $\lambda$ and $Z_{\mathbb{\bC}}(\lambda):=\mathfrak{A}_{\bC} \bf{v}^+$ the cyclic $\mathfrak{A}_{\mathbb{C}}$-module associated to $\lambda$ (called {\em Verma module} in Humphreys \cite{Humphreys2}). We recall that ${\bf v}^+$ is uniquely determined modulo scalars by the following conditions: $H{\bf v}^+ = \lambda(H){\bf v}^+$ for every $H\in \mathfrak{H}$ and $X_\alpha{\bf v}^+ = 0$ for every $\alpha\in \Phi^+$.

We have $Z_{\mathbb{\bC}}(\lambda)=\oplus_{\mu\in\Lambda,~\mu\preceq \lambda} {\bf V}_{\mu}$ where ${\bf V}_{\mu}$ is the {\em weight space} with weight $\mu$ (Humphreys~\cite[Chapter 20]{Humphreys}). More explicitly, $\mu = \lambda - \sum_{i=1}^nm_i\alpha_i$ for nonnegative integers $m_1,..., m_n$ and ${\bf V}_\mu$ is the direct sum of the $1$-dimensional subspaces $\langle A^-_t{\bf v}^+\rangle$ for $A^-_t$ such that $d_i(A^-_t) = -m_i$ for every $i = 1, 2,..., n$. Moreover, ${\bf V}_\mu$ is the $\mu$-eigenspace of $\mathfrak{H}$, namely its vectors are the vectors ${\bf x}\in Z_{\mathbb{C}}(\lambda)$ such that $H({\bf x}) = \mu(H){\bf x}$ for every $H\in \mathfrak{H}$. In particular, ${\bf V}_\lambda = \langle {\bf v}^+\rangle$. Note also that, for every weight $\mu \preceq \lambda$ and every  element $X \in{\cal A}$ we have $X{\bf V}_\mu \subseteq {\bf V}_{\mu + \sum_{i=1}^nd_i(X)\alpha_i}$, with the convention that ${\bf V}_{\mu + \sum_{i=1}^nd_i(X)\alpha_i}= 0$ if $\lambda \prec\mu + \sum_{i=1}^nd_i(X)\alpha_i$.

The module $Z_{\bC}(\lambda)$ admits a unique maximal proper submodule $J_{\bC}(\lambda)$. The {\em Weyl module} relative to $\lambda$ is defined as the quotient $V_{\bC}(\lambda):= Z_{\bC}(\lambda)/J_{\bC}(\lambda)$. The Weyl module $V_{\bC}(\lambda)$ is finite dimensional. In particular, $J_{\bC}(\lambda)$ contains all weight spaces ${\bf V}_\mu$ with $\mu \prec 0$. It follows that for every $\alpha\in \Phi^+$ both elements $X_\alpha$ and $Y_\alpha$ act on $V_{\bC}(\lambda)$ (by left multiplication) as nilpotent mappings. Thus, the exponential map $e^{X_\alpha t} := \sum_{i=0}^\infty\frac{X_\alpha^it^i}{i!}$ can be defined on $V_{\bC}(\lambda)$ for every $\alpha\in \Phi$ (and every $t\in\bC$). This map is (linear and) invertible. If ${\bf X}$ denotes the Dynkin type of the root system $\Phi$, the group $G_{\bC}(\lambda) = \langle e^{X_\alpha t}\rangle_{\alpha\in\Phi, t\in \bC}$ is the complex Chevalley group of type ${\bf X}$ associated to the weight $\lambda$.

Before to go on we fix a notation which allows us to immediately recognize if an element is regarded as an element of $V_{\bC}(\lambda)$ or $Z_{\bC}(\lambda)$.  We denote vectors and subspaces of $Z_{\mathbb{\bC}}(\lambda)$ by bold characters, as we have done so far, keeping normal font letters for the corresponding vectors and subspaces of $V_{\bC}(\lambda)$. Thus, a letter as $\bf{v}$ stands for a vector of $Z_{\mathbb{\bC}}(\lambda)$ while $v$ denotes the vector $\pi_{\bC}({\bf v})$, where $\pi_{\bC}$ is the canonical projection of $Z_{\bC}(\lambda)$ onto $V_{\bC}(\lambda)$. Similarly, if ${\bf V}_\mu$ is a weight space of $Z_{\bC}(\lambda)$ then $V_\mu = \pi_\bC({\bf V}_\mu)$.

We shall now define the analogues $G_{\F}(\lambda)$ and $V_\F(\lambda)$ of $G_{\bC}(\lambda)$ and $V_{\bC}(\lambda)$  over an arbitrary field $\F$. To that goal we must preliminary define suitable $\mathbb{Z}$-analogues $Z_{\mathbb{\mathbb{Z}}}(\lambda)$, $J_{\mathbb{Z}}(\lambda)$ and $V_{\mathbb{Z}}(\lambda)$ of $Z_{\mathbb{\bC}}(\lambda)$, $J_{\bC}(\lambda)$ and $V_{\bC}(\lambda)$.

The definition of $Z_{\mathbb{Z}}(\lambda)$ is obvious: we put $Z_{\mathbb{Z}}(\lambda) := \mathfrak{A}_{\mathbb{Z}}{\bf v}^+$ ($= \mathfrak{A}^-_{\mathbb{Z}}{\bf v}^+$ in view of the previous description of $Z_{\bC}(\lambda)$ and its weight spaces). In order to define $J_{\mathbb{Z}}(\lambda)$ we must firstly define $V_{\mathbb{Z}}(\lambda)$.

According to Steinberg~\cite[page 17]{Stein}, there exists a basis $B$ of  $V_{\bC}(\lambda)$ such that $\langle B \rangle_{\bZ}=\mathfrak{A}^-_{\mathbb{Z}} v^+$. The equality $\langle B\rangle_{\mathbb{Z}} = \mathfrak{A}^-_{\mathbb{Z}}v^+$ implies that we can choose the basis $B$ of $\mathfrak{A}^-_{\mathbb{Z}}v^+$ in such a way that its elements are monomial vectors of the form $A^-_tv^+$, namely $B\subseteq {\cal A}^-v^+$. So, we assume that $B$ has been chosen in that way. We put $V_{\mathbb{Z}}(\lambda) := \langle B\rangle_{\mathbb{Z}} = \mathfrak{A}^-_{\mathbb{Z}}v^+$ ($= \mathfrak{A}_{\mathbb{Z}}v^+$, since $\mathfrak{A}^-_{\mathbb{Z}}{\bf v}^+ = \mathfrak{A}_{\mathbb{Z}}{\bf v}^+$).

We have $V_{\mathbb{Z}}(\lambda) =\oplus_{\mu\preceq \lambda} V_{\mu,\mathbb{Z}}$, where $V_{\mu,\mathbb{Z}} = \pi_{\bC}({\bf V}_{\mu,\mathbb{Z}})$ and ${\bf V}_{\mu,\mathbb{Z}} = {\bf V}_\mu\cap Z_{\mathbb{Z}}(\lambda)$. Denoted by $\overline{\bf B}_\mu$ the monomial basis of ${\bf V}_{\mu}$ contained in ${\cal A}^-{\bf v}^+$, the set $\overline{\bf B}_{\mu}$ contains a subset ${\bf B}_{\mu}$ such that $\pi_\bC$ induces a bijection from ${\bf B}_{\mu}$ to a basis $B_\mu = \pi_{\bC}(\B_\mu)$ of $V_\mu$ contained in $B$. (We warn that ${\bf B}_\mu = \emptyset$ is allowed; this is certainly the case when $\mu \prec 0$, since $V_\mu = 0$ in this case.)

We shall now replace the set $\overline{\B}_\mu\setminus \B_{\mu}$ with an independent set of vectors $\B'_\mu\subset J_{\mathbb{C}}(\lambda)\cap\langle \overline{\B}_\mu\rangle_{\mathbb{Z}}$ (where $\langle\overline{\B}_\mu\rangle_{\mathbb{Z}}$ stands for the set of integral combinations of vectors of $\overline{\B}_\mu$) such that $\langle\B_{\mu}\cup \B'_{\mu}\rangle_{\mathbb{Z}} = \langle \overline{\B}_\mu\rangle_{\mathbb{Z}}$ (whence $\B_{\mu}\cup\B'_{\mu}$ is a basis of ${\bf V}_{\mu}$).  Having done this we will be able to define $J_{\mathbb{Z}}(\lambda)$.

Let ${\bf b}'_1, {\bf b}'_2,..., {\bf b}'_r$ be the elements of $\overline{\B}_\mu\setminus\B_\mu$, taken in some given order. We construct $\B'_{\mu}$ recursively as follows. We put $\B'_{\mu, 0} := \emptyset$. Assuming to have already defined $\B'_{\mu,k}$ for $0\leq k < r$, if $\pi_\bC({\bf b}'_{k+1}) = 0$ then $\B'_{\mu,k+1}: = \B'_{\mu, k}\cup\{{\bf b}'_{k+1}\}$. Otherwise $b'_{k+1} := \pi_\bC({\bf b}_{k+1})\in V_\mu\cap{\cal A}^-v^+ \subset V_{\mu, \mathbb{Z}}$ is an integral combination of the vectors $b_1, b_2,..$ of $B_\mu$, say $b'_{k+1} = \sum_iz_ib_i$. Clearly ${\bf b}'_{k+1}-\sum_iz_i{\bf b}_i \in J_{\bC}(\lambda)$. We put $\B'_{\mu,k+1}: = \B'_{\mu, k}\cup\{{\bf b}'_{k+1}-\sum_iz_i{\bf b}_i\}$. Finally, $\B'_\mu := \B'_{\mu, r}$.

Put $J_{\mathbb{Z}}(\lambda):=\langle \cup_{\mu\preceq\lambda} \B'_{\mu}\rangle_{\bZ}$. Then $(\cup_{\mu\preceq\lambda}\B_{\mu})\cup(\cup_{\mu\preceq\lambda}\B'_{\mu})$ is a basis of $Z_{\mathbb{\bZ}}(\lambda)$ and $J_{\mathbb{Z}}(\lambda) = Z_{\mathbb{Z}}(\lambda)\cap J_\bC(\lambda)$ (as it follows from the definition of $\B_\mu$ and the properties of the sets $\B'_\mu$ defined above).  We can now interpret $V_{\mathbb{Z}}(\lambda)$ (previously defined as $V_{\mathbb{Z}}(\lambda):=\langle B\rangle_{\bZ} = \mathfrak{A}_{\mathbb{Z}}v^+ = \mathfrak{A}^-_{\mathbb{Z}}v^+$) as the quotient $V_{\mathbb{Z}}(\lambda):=Z_{\bZ} (\lambda)/J_{\mathbb{Z}}(\lambda)$.

The projection $\pi_\bC:Z_{\bC}(\lambda)\rightarrow V_{\bC}(\lambda)$ induces a homomorphism $\pi_{\mathbb{Z}}$ from $Z_{\mathbb{Z}}(\lambda)$ to $V(\lambda)_{\mathbb{Z}}$ and $J_{\mathbb{Z}}(\lambda) = \mathrm{ker}(\pi_{\mathbb{Z}})$.

Let now $\F$ be an arbitrary field. Then we put $\mathfrak{A}_{\F}:=\F\otimes_{\mathbb{Z}} \mathfrak{A}_{\bZ}$, $\mathfrak{A}^-_\F := \F\otimes_{\mathbb{Z}}\mathfrak{A}^-_{\mathbb{Z}}$, $Z_\F(\lambda) = \F\otimes_\bZ Z_\bZ(\lambda)$,  $V_\F(\lambda) = \F\otimes_\bZ V_{\mathbb{Z}}(\lambda)$ and $J_\F(\lambda) = \F\otimes_\bZ J_\bZ(\lambda)$.  Since $V_\bZ(\lambda) = Z_\bZ(\lambda)/J_\bZ(\lambda)$, we also have $V_\F(\lambda) = Z_{\mathbb{F}}(\lambda)/J_{\mathbb{F}}(\lambda)$. The Chevalley group $G_\F(\lambda)$ is
defined as a group of linear transformations of $V_\F(\lambda)$,  by means of the exponentials $e^{X_\alpha t}$ just as in the complex case. The Lie algebra $\mathfrak{L}_\F = \F\otimes_\mathbb{Z}\mathfrak{L}_\mathbb{Z}$ is the Lie algebra of $G_\F(\lambda)$.

We finish this subsection with a result on homomorphic images of Weyl modules, already but implicitly used in Section \ref{Introduction}.

\begin{theo}\label{immagini}
Let $U$ be a finite dimensional  $G_\F(\lambda)$-module over $\F$, satisfying the following properties.

\begin{itemize}
\item[(1)] $U$ is also an $\mathfrak{A}_\F$-module and the action of $G_\F(\lambda)$ on $U$ is induced by the action of $\mathfrak{A}_\F$. Explicitly, for every root $\alpha\in\Phi$ the element $X_\alpha$ acts on $U$ as a nilpotent endomorphism and the element $x_\alpha(t) = e^{X_\alpha t}$ of $G_\F(\lambda)$ maps $v$ onto $\sum_{n=0}^\infty t^n\frac{X_{\alpha}^n}{n!}(v)$, for every vector $v\in U$.
\item[(2)] $U$ admits a highest weight vector of weight $\lambda$, namely a vector $\bar{v}^+$ such that $X_\alpha(\bar{v}^+) = 0$ for every $\alpha\in\Phi^+$ and $H_i(\bar{v}^+) = \lambda(H_i){\bar{v}^+}$ for every $i = 1, 2,..., n$.
\end{itemize}
Then there exists a unique homomorphism of $\mathfrak{A}_\F$-modules from $V_\F(\lambda)$ to $U$ mapping $v^+$ onto $\bar{v}^+$ (whence $g(v^+)$ onto $g(\bar{v}^+)$ for every element $g\in\mathfrak{A}_\F$).
\end{theo}

This theorem is implicit in the final result of Humphreys \cite[3.3]{Humphreys2}. In \cite[3.3]{Humphreys2} the field $\F$ is assumed to be algebraically closed, but in the hypotheses of Theorem \ref{immagini} that assumption can be dropped. Indeed we can firstly replace $\F$ with its algebraic closure $\overline{\F}$ and apply the theorem of \cite[3.3]{Humphreys2} to $G_{\overline{\F}}(\lambda)$ acting on $V_{\overline{\F}}(\lambda)$ and $\overline{\F}\otimes_\F U$. Next we turn back to $G_\F(\lambda)$ acting on $V_\F(\lambda)$ and $U$.

\subsection{A Chevalley basis of $\mathfrak{o}(2n+1,\F)$}\label{orthogonal algebra}

Given a field $\F$, let $\mathfrak{L}^\circ_{\mathbb{F}} = \mathfrak{o}(2n+1,\mathbb{F}) = \mathfrak{H}\oplus(\oplus_{\alpha\in \Phi}\mathfrak{X}_{\alpha})$ be the Lie algebra of type $B_n$ over $\F$, where $\mathfrak{H}$, $\Phi$ and $\mathfrak{X}_{\alpha}$ have the same meaning as in Section \ref{general setting}. Let $\mathfrak{A}^\circ_{\F}$ be the enveloping algebra of $\mathfrak{L}^\circ_{\F}$, $\Pi = \{\alpha_1,..., \alpha_{n}\}$ a given basis of simple roots of $\Phi$ and $\{X_\alpha, Y_\alpha\}_{\alpha\in \Phi^+}\cup\{H_i\}_{i=1}^n$ a Chevalley basis of $\mathfrak{L}^\circ_{\F}$ where $Y_\alpha := X_{-\alpha}$, as in Section \ref{general setting}.

The roots $\alpha\in \Phi$ are represented by the following vectors of $\mathbb{R}^n$, where $(u_i)_{i=1}^{n}$ is the natural orthonormal basis of $\mathbb{R}^n$ (see Humphreys~\cite[12.1]{Humphreys}):
\[(\Pi) ~
\left\{\begin{array}{ll}
\alpha_i=u_i-u_{i+1}~ \mbox{for}~1\leq i\leq n-1 &\mbox{(long simple roots)}\\
\alpha_n=u_n & \mbox{(short simple root)}
\end{array}\right.\]
\[(\Phi^+) ~ \left\{\begin{array}{ll}
u_i-u_j = \sum_{r=i}^{j-1}\alpha_r ~ \mbox{for}~1\leq i< j\leq n &\mbox{(long)}\\
{} & \\
u_i = \sum_{r=i}^n\alpha_r ~ \mbox{for}~ i = 1, 2,..., n & \mbox{(short)}\\
{} & \\
u_i+u_j = \sum_{r=i}^{j-1}\alpha_r + 2\sum_{r=j}^n\alpha_r ~\mbox{for}~1\leq i<j\leq n & \mbox{(long)}
\end{array}\right.\]
\[(\Phi^-) ~ \left\{\begin{array}{ll}
u_j-u_i = -\sum_{r=i}^{j-1}\alpha_r ~ \mbox{for}~1\leq i< j\leq n &\mbox{(long)}\\
{} & \\
-u_i = -\sum_{r=i}^n\alpha_r ~ \mbox{for}~ i = 1, 2,..., n & \mbox{(short)}\\
{} & \\
-u_i-u_j = -\sum_{r=i}^{j-1}\alpha_r - 2\sum_{r=j}^n\alpha_r ~\mbox{for}~1\leq i<j\leq n & \mbox{(long)}
\end{array}\right.\]
As in Section \ref{Introduction}, let $V = V(2n+1,\F)$ and let $\eta$ be a non-singular quadratic form on $V$. Up to rescaling the form $\eta$, we can assume to have chosen an ordered basis $E = (e_1,e_2,\dots,e_{2n+1})$ of $V$ with respect to which
\begin{equation}\label{eta}
\eta(x_1,..., x_{2n+1}) = \sum_{i=1}^nx_ix_{n+i} + x_{2n+1}^2.
\end{equation}
Regarding $G = \mathrm{SO}(2n+1,\F)$ as the subgroup of $\mathrm{SL}(V)$ preserving $\eta$ and $\mathfrak{L}_\F^\circ$ as the Lie algebra of $G$, the algebra $\mathfrak{L}^\circ_\F$ acts on $V$ as a subalgebra of $\mathrm{End}(V)$.

We can assume that the chosen Chevalley basis $\{X_\alpha, Y_\alpha\}_{\alpha\in \Phi^+}\cup\{H_i\}_{i=1}^n$ of $\mathfrak{L}^\circ
_\F$ acts as follows on the elements of $E$, where $\delta_{i,j}$ is the Kronecker symbol:
\[\begin{array}{rcll}
X_{u_i-u_j}(e_l) & = & \delta_{l,j}e_i - \delta_{l,n+i}e_{n+j}, & (1\leq i < j \leq n),\\
Y_{u_i-u_j}(e_l) & = & \delta_{l,i}e_j - \delta_{l,n+j}e_{n+i}, & (1\leq i < j \leq n),\\
X_{u_i+u_j}(e_l) & = & \delta_{l,n+i}e_j- \delta_{l,n+j}e_i, & (1\leq i < j \leq n),\\
Y_{u_i+u_j}(e_l) & = & -\delta_{l,i}e_{n+j} + \delta_{l,j}e_{n+i}, & (1\leq i < j \leq n),\\
X_{u_i}(e_l) & = & \delta_{l,n+i}e_{2n+1} - 2\delta_{l,2n+1}e_i, & (1\leq i \leq n),\\
Y_{u_i}(e_l) & = & -\delta_{l,i}e_{2n+1} + 2\delta_{l,2n+1}e_{n+i}, & (1\leq i \leq n),
\end{array}\]
\[\begin{array}{ll}
H_i(e_l)  =  \delta_{l,i}e_i - \delta_{l,i+1}e_{i+1} - \delta_{l,n+i}e_{n+i}+\delta_{l,n+i+1}e_{n+i+1}, & (1\leq i < n),\\
H_n(e_l) = 2\delta_{l,n}e_n - 2\delta_{l,2n}e_{2n}.
\end{array}\]
(Recall that $H_i = H_{\alpha_i} = H_{u_i-u_{i+1}}$ if $i < n$ and $H_n = H_{\alpha_n} = H_{u_n}$.) Moreover,
\[\begin{array}{ll}
H_{u_i-u_j} = \sum_{r=i}^{j-1}H_r = -H_{u_j-u_i}, & (1\leq i < j\leq n),\\
{} & \\
H_{u_i} = 2\sum_{r=i}^{n-1}H_r + H_n = -H_{-u_i}, & (1\leq i \leq n)\\
{} & \\
H_{u_i+u_j} = \sum_{r=1}^{j-1}H_r + 2\sum_{r=j}^{n-1}H_r + H_n = -H_{-u_i-u_j}, & (1\leq i<j\leq n).
\end{array}\]
The elements $X_\alpha$ and $Y_\alpha$ act on $V$ as nilpotent endomorphisms of exponent 2 or 3, according to whether $\alpha$ is long or short. Moreover,
\[\begin{array}{rcll}
\frac{X^2_{u_i}}{2}(e_l) & = & -\delta_{l,n+i}e_{i}  & (1\leq i \leq n),\\
{} & & & \\
\frac{Y^2_{u_i}}{2}(e_l) & = & -\delta_{l,i}e_{n+i}  & (1\leq i \leq n).
\end{array}\]

We now turn to $W_{k} = \bigwedge^kV$. We firstly fix some notation. We set $I :=  \{1, 2,..., 2n+1\}$ and $E_\wedge := (e_J)_{J\in{I\choose k}}$, where $I\choose k$ stands for the
set of subsets of $I$ of size $k$ and $e_J = e_{j_1}\wedge e_{j_2}\wedge...\wedge e_{j_k}$ for every $k$-subset $J = \{j_1,..., j_k\}$ of $I$, with the convention that $j_1 < j_2 < ... < j_k$. Thus, $E_\wedge$ is the basis of $W_{k}$ associated to $E$.

The action of $\mathfrak{L}^\circ_\F$ on $V$ induces an action of $\mathfrak{L}^\circ_\F$ on $W_{k}$ defined as follows: if $g\in \mathfrak{L}^\circ_\F$ and $v = v_1\wedge v_2\wedge...\wedge v_k$ then
\begin{equation}\label{azione su wedge}
g(v) = \sum_{i=1}^kv_1\wedge...\wedge v_{i-1}\wedge g(v_i)\wedge v_{i+1}\wedge... \wedge v_k.
\end{equation}
The action of the elements $X_\alpha^t/t!$ and $Y_\alpha^t/t!$ on $W_{k}$ can be computed with the help of formula (\ref{azione su wedge}). We leave these computations for the reader.

Let $\bar{v}^+:=e_1\wedge e_2\wedge\dots\wedge e_k$. Then $\bar{v}^+$ behaves as the maximal vector $v^+$ of $V(\lambda_k^\circ)$ (where $\lambda_k^\circ$ has the meaning stated in Section
\ref{Introduction}), namely $X_{\alpha}(\bar{v}^+)=0$ for any $\alpha\in \Phi^+$ and $H(\bar{v}^+)
= \lambda_k^\circ(H)\bar{v}^+$ for every $H\in \mathfrak{H}$. In particular,
\[\begin{array}{ccc}
H_{u_i-u_{i+1}}(\bar{v}^+) = \left\{\begin{array}{l}
\bar{v}^+ ~ \mbox{if $i = k$},\\
0 ~ \mbox{otherwise}
\end{array}\right. & &
H_{u_n}(\bar{v}^+) = \left\{\begin{array}{l}
2\bar{v}^+ ~ \mbox{if $k = n$},\\
0 ~ \mbox{otherwise}
\end{array}\right.
\end{array}\]
Note that
\begin{equation}\label{lambda kappa}
\lambda_k^\circ = \sum_{i=1}^{k-1}i\alpha_i + k\cdot\sum_{i=k}^n\alpha_i.
\end{equation}
(Humphreys \cite[13.2]{Humphreys}). Thus,
\[\lambda^\circ_k(H_i) = \langle \lambda_k^\circ, \alpha_i\rangle = \left\{\begin{array}{ll}
-(i-1)+2i-(i+1) = 0 & \mbox{if $i < k < n$},\\
-(k-1) + 2k - k = 1 & \mbox{if $i = k < n$},\\
-k+2k-k = 0 & \mbox{if $k < i < n$},\\
2k-2k = 0 & \mbox{if $k < i = n$},\\
-2(n-1)+ 2n = 2 & \mbox{if $k = i = n$},
\end{array}\right.\]
which fits with the above values of $H_{u_i-u_{i+1}}(\bar{v}^+)$ and $H_{u_n}(\bar{v}^+)$ (recall that $H_i = H_{u_i-u_{i+1}}$ if $i < n$ and $H_n = H_{u_n}$).

Let $W^\circ_{k}(\F) := \mathfrak{A}^\circ_{\F}(\bar{v}^+)$. Thus $W^\circ_{k}(\F)$ is the module called $W^\circ_{k}$ in Section \ref{Introduction}, but now we keep a record of the field $\F$ in our notation. In view of the above, the homomorphism $\varphi_{k,\F}:V_\F(\lambda^\circ_k)\rightarrow W^\circ_{k}(\F)$ (Theorem \ref{immagini}) can be chosen in such a way that $\varphi_{k,\F}(v^+) = \bar{v}^+$. This condition uniquely determines $\varphi_{k,\F}$.

As recalled in Section \ref{Introduction}, the homomorphism $\varphi_{k,\F}$ is an isomorphism if and only if either $\mathrm{char}(\F) \neq 2$ or $k = 1$. If $\mathrm{char}(\F) = 2$ and $k > 1$ then the kernel $K_k$ of $\varphi_{k.\F}$ has dimension $\mathrm{dim}(K_k) = {{2n+1}\choose{k-2}}$.

\subsection{A Chevalley basis of $ \mathfrak{sp}(2n,\mathbb{\F})$}\label{symplectic algebra}

Let $\mathfrak{L}^{\mathrm{sp}}_{\mathbb{F}} = \mathfrak{sp}(2n,\mathbb{F}) = {\mathfrak{C}}\oplus(\oplus_{\beta\in\overline{\Phi}}{\mathfrak{U}}_{\beta})$ be the Lie algebra of type $C_n$, where ${\mathfrak{C}}$ is the Cartan subalgebra, $\overline{\Phi}$ is the root system of type $C_n$ and $\mathfrak{U}_{\beta}$ is the $1$-dimensional subalgebra of $\mathfrak{L}^{\mathrm{sp}}_{\mathbb{F}}$ corresponding to the root $\beta$. Let $\mathfrak{A}^{\mathrm{sp}}_{\bC}$ be the enveloping algebra of $\mathfrak{L}^{\mathrm{sp}}_{\bC}$. Chosen a basis $\overline{\Pi} = \{\beta_1,..., \beta_{n}\}$ of simple roots for $\overline{\Phi}$ let $\overline{\Phi}^+$ (respectively $\overline{\Phi}^-$) be the set of roots that are positive (negative) with respect to $\overline{\Pi}$. Let $\{C_i\}_{i=1}^n\cup \{U_\beta, V_\beta\}_{\beta\in \overline{\Phi}}$ be a Chevalley basis of $\mathfrak{L}^{\mathrm{sp}}_{\F}$, where $C_i\in \mathfrak{C}$ for $i = 1,...,n$ and  $U_\beta\in\mathfrak{U}_{\beta}$ and $V_\beta = U_{-\beta}\in\mathfrak{U}_{-\beta}$ for $\beta\in \overline{\Phi}^+$.

The positive roots $\beta\in \overline{\Phi}^+$ are represented by the following vectors of $\mathbb{R}^n$ where $(u_i)_{i=1}^{n}$ is the natural orthonormal basis of $\mathbb{R}^n$, as in the previous subsection (see Humphreys~\cite[12.1]{Humphreys}):
\[\begin{array}{ll}
\beta_i=u_i-u_{i+1} ~ \mbox{for}~1\leq i\leq n-1 &\mbox{(short simple)}\\
{} & \\
\beta_n= 2u_n & \mbox{(long simple)} \\
{} & \\
u_i-u_j = \sum_{r=i}^{j-1}\beta_r ~ \mbox{for}~1\leq i< j\leq n &\mbox{(short)}\\
{} & \\
2u_i = 2\sum_{r=i}^{n-1}\beta_r + \beta_n ~ \mbox{for}~ i = 1, 2,..., n  & \mbox{(long)}\\
{} & \\
u_i+u_j = \sum_{r=i}^{j-1}\beta_r + 2\sum_{r=j}^{n-1}\beta_r + \beta_n ~\mbox{for}~1\leq i<j\leq n & \mbox{(short)}
\end{array}\]
The negative roots are the opposite of these. We omit to list them. Note that the short roots of $\overline{\Phi}$ are just the long roots of the root system $\Phi$ of type $B_n$ while the long roots of $\overline{\Phi}$ are the short roots of $\Phi$ multiplied by $2$.

With $V = V(2n+1,\F)$ and $E = \{e_1,..., e_{2n+1}\}$ as in the previous subsection, let $\overline{V}$ be the hyperplane of $V$ spanned by the subset $\overline{E}=\{e_1,e_2,\dots,e_{2n}\}$ of $E$. Let $\overline{V}$ be endowed with the non-degenerate alternating form $\varsigma$ represented as follows with respect to the basis $\overline{E}$ of $\overline{V}$:
\begin{equation}\label{varsigma}
\varsigma((x_1,x_2,\dots, x_{2n}), (y_1,y_2,\dots, y_{2n})) = \sum_{i=1}^{n}x_iy_{n+i} -\sum_{i=1}^{n}x_{n+i}y_i.
\end{equation}
(Note that if $\mathrm{char}(\F) = 2$ then, according to (\ref{eta}), $\varsigma$ is the form induced on $\overline{V}$ by the bilinearization of the quadratic form $\eta$ of $V$.) Regarded $\overline{G} := \mathrm{Sp}(2n,\F)$ as the subgroup of $\mathrm{SL}(\overline{V})$ preserving $\varsigma$ and $\mathfrak{L}_\F^{\mathrm{sp}}$ as the Lie algebra of $\overline{G}$, the algebra $\mathfrak{L}^{\mathrm{sp}}_\F$ acts on $\overline{V}$ as a subalgebra of $\mathrm{End}(\overline{V})$.

We can assume that the chosen Chevalley basis $\{U_\beta, V_\beta\}_{\beta\in \overline{\Phi}^+}\cup\{C_i\}_{i=1}^n$ of $\mathfrak{L}^{\mathrm{sp}}_\F$ acts as follows on the elements of $\overline{E}$:
\[\begin{array}{rcll}
U_{u_i-u_j}(e_l) & = & \delta_{l,j}e_i - \delta_{l,n+i}e_{n+j}, & (1\leq i < j \leq n),\\
V_{u_i-u_j}(e_l) & = & \delta_{l,i}e_j - \delta_{l,n+j}e_{n+i}, & (1\leq i < j \leq n),\\
U_{u_i+u_j}(e_l) & = & \delta_{l,n+i}e_j+ \delta_{l,n+j}e_i, & (1\leq i < j \leq n),\\
V_{u_i+u_j}(e_l) & = & \delta_{l,i}e_{n+j} + \delta_{l,j}e_{n+i}, & (1\leq i < j \leq n),\\
U_{2u_i}(e_l) & = & \delta_{l,n+i}e_{i}, & (1\leq i \leq n),\\
V_{2u_i}(e_l) & = & \delta_{l,i}e_{n+i}, & (1\leq i \leq n),
\end{array}\]
\[\begin{array}{ll}
C_i(e_l)  =  \delta_{l,i}e_i - \delta_{l,i+1}e_{i+1} - \delta_{l,n+i}e_{n+i}+\delta_{l,n+i+1}e_{n+i+1}, & (1\leq i < n),\\
C_n(e_l) = \delta_{l,n}e_l - \delta_{l,2n}e_{2n}.
\end{array}\]
Note that $U_\beta$ and $V_\beta$ act on $\overline{V}$ as nilpotent endomorphisms of exponent 2 for every $\beta\in \overline{\Phi}^+$.  Moreover, the elements $X_{u_i-u_j}$, $Y_{u_i-u_j}$, $\frac{X_{u_i}^2}{2}$, $\frac{Y_{u_i}^2}{2}$, $X_{u_i+u_j}$ and $Y_{u_i+u_j}$ of $\mathfrak{A}^\circ_\F$ stabilize the subspace $\overline{V}$ of $V$  and the following hold for every vector $v\in \overline{V}$
\begin{equation}\label{corrispondenza1}
\left.\begin{array}{ll}
U_{u_i-u_j}(v) = X_{u_i-u_j}(v), &  V_{u_i-u_j}(v) = Y_{u_i-u_j}(v),\\
{} & \\
U_{2u_i}(v) = -\frac{X^2_{u_i}}{2}(v), & V_{2u_i}(v) = -\frac{Y_{u_i}^2}{2}(v),
\end{array}\right\}
\end{equation}
\begin{equation}\label{corrispondenza2}
\left.\begin{array}{lcl}
U_{u_i+u_j}(v) & = & X_{u_i+u_j}(v) -2X_{e_i-e_j}\frac{X_{u_j}^2}{2}(v) = \\
 {} & = & X_{u_i+u_j}(v) + 2U_{e_i-e_j}U_{2u_j}(v),\\
{} & & \\
V_{u_i+u_j}(v) & = & Y_{u_i+u_j}(v) -2\frac{Y_{u_j}^2}{2}Y_{e_i-e_j}(v) = \\
 {} & = & Y_{u_i+u_j}(v)  +2V_{2u_j}V_{e_i-e_j}.
\end{array}\right\}
\end{equation}
In particular, when $\mathrm{char}(\F) = 2$
\begin{equation}\label{corrispondenza3}
\begin{array}{ll}
U_{u_i+u_j}(v) = X_{u_i+u_j}(v), & V_{u_i+u_j}(v) = Y_{u_i+u_j}(v).
\end{array}
\end{equation}
Moreover, $U_{u_i+u_j} = [U_{u_i-u_j}, U_{2u_j}]$ and $V_{u_i+u_j} = [V_{2u_j}, V_{u_i-u_j}]$. Consequently,
\begin{equation}\label{commutazione}
\left.\begin{array}{l}
U_{u_i+u_j}(v) =  \frac{X^2_{u_j}}{2}X_{u_i-u_j}(v) - X_{u_i-u_j}\frac{X^2_{u_j}}{2}(v),\\
{}\\
V_{u_i+u_j} = Y_{u_i-u_j}\frac{Y_{u_j}^2}{2}(v) - \frac{Y_{u_j}^2}{2}Y_{u_i-u_j}(v).
\end{array}\right\}
\end{equation}

We now turn to $\overline{W}_{k} = \bigwedge^k\overline{V}$. We take $\overline{E}_\wedge := (e_J)_{J\in{\bar{I}\choose k}}$ as a basis for $\overline{W}_{k}$, where $\bar{I} := \{1,2,..., 2n\}$. The action of $\mathfrak{L}^{\mathrm{sp}}_\F$ on $\overline{V}$ induces an action of $\mathfrak{L}^{\mathrm{sp}}_\F$ on $\overline{W}_{k}$, defined according to the rule (\ref{azione su wedge}).

Note firstly that $U_\beta$ and $V_\beta$ act on $\overline{W}_{k}$ as nilpotent endomorphisms of exponent 3 or 2 according to whether $\beta$ is short or long. In view of (\ref{corrispondenza1}) the elements $U_{u_i-u_j}$, $\frac{U_{u_i-u_j}^2}{2}$, $V_{u_i-u_j}$, $\frac{V_{u_i-u_j}^2}{2}$, $U_{2u_j}$ and $V_{2u_j}$ act on $\overline{W}_{k}$ as $X_{u_i-u_j}$, $\frac{X_{u_i-u_j}^2}{2}$, $Y_{u_i-u_j}$, $\frac{Y_{u_i-u_j}^2}{2}$, $-\frac{X_{u_i}^2}{2}$ and $-\frac{Y_{u_i}^2}{2}$ respectively. According to (\ref{corrispondenza3}), denoted by $\langle \overline{E}_\wedge\rangle_{2\mathbb{Z}}$ the set of integral combination of vectors of $\overline{E}_\wedge$ with even coefficients,
\begin{equation}\label{corrispondenza4}
\{U_{u_i+u_j}(e_J)-X_{u_i+u_j}(e_J),~ V_{u_i+u_j}(e_J)-Y_{u_i+u_j}(e_J)\} ~\subset ~\langle \overline{E}_\wedge\rangle_{2\mathbb{Z}}
\end{equation}
for every $e_J\in \overline{E}_\wedge$. In particular, if $\mathrm{char}(\F) = 2$ then $U_{u_i+u_j}$ and $V_{u_i+u_j}$ act on $\overline{W}_{k}$ in the same way as $X_{u_i+u_j}$ and $Y_{u_i+u_j}$ respectively.

Let $\bar{v}^+:=e_1\wedge e_2\wedge\dots\wedge e_k$, as in the previous subsection and let $\overline{W}^{\mathrm{sp}}_{k}(\F) = \mathfrak{A}_\F^{\mathrm{sp}}(\bar{v}^+)$. Thus,  $\overline{W}^{\mathrm{sp}}_{k}(\F)$ is the module called $\overline{W}_{k}^{\mathrm{sp}}$ in Section \ref{Introduction}. By (\ref{corrispondenza1}), (\ref{corrispondenza2}) and (\ref{corrispondenza3}), we see that $\overline{W}^{\mathrm{sp}}_{k}(\F)\subseteq W_{k}^\circ(\F)$.

Let $\lambda_1^{\mathrm{sp}}, \lambda_2^{\mathrm{sp}},..., \lambda_n^{\mathrm{sp}}$ be the fundamental  dominant weights for the root system of type $C_n$. It is well known that $\overline{W}^{\mathrm{sp}}_{k}(\F)$ is isomorphic to the Weyl module $V_\F(\lambda^{\mathrm{sp}}_k)$ and that, if $v^+$ is the highest weight vector chosen in $V_\F(\lambda^{\mathrm{sp}}_k)$, we can assume that the isomorphism from $V_\F(\lambda^{\mathrm{sp}}_k)$ to $\overline{W}_{k}^{\mathrm{sp}}$ maps $v^+$ onto $\bar{v}^+$ (see e.g. Premet and Suprunenko \cite{PS}).

\section{Proof of Theorem \ref{main theorem}}\label{section-main theorem}

Throuhgout this section $\mathrm{char}(\F) = 2$. If $k = 1$ then $K_k = 0$ and ${\cal N}_k$ is 1-dimensional. In fact ${\cal N}_1$ is the nucleus of the quadratic form $\eta$ chosen to define $G$. Thus, in this case the first claim of Theorem \ref{main theorem} is obvious. The second claim is empty.

Assume $k > 1$. We keep the notation of Subsections \ref{orthogonal algebra} and \ref{symplectic algebra}. In particular, if $v^+$ is the highest weight vector chosen in $V_\F(\lambda^\circ_k)$ and $\bar{v}^+ = e_1\wedge...\wedge e_k$, we can assume to have chosen the homomorphism $\varphi_{k,\F}:V_\F(\lambda^\circ_k)\rightarrow W_{k}^\circ(\F)$ in such a way that it maps $v^+$ onto $\bar{v}^+$. We put $K_k := \mathrm{ker}(\varphi_{k,\F})$. Clearly,  $K_k$ is an $\mathfrak{A}^\circ_\F$-submodule of $V_\F(\lambda^\circ_k)$, namely

\begin{lemma}\label{lemma0}
We have $\mathfrak{A}^\circ_\F(K_k) = K_k$.
\end{lemma}

Let $v_1^+ : = Y_{u_k}(v^+)$, $V_1(\F) := \mathfrak{A}^\circ_\F(v_1^+)$ and
\[\bar{v}_1^+ := \varphi_{k,\F}(v_1^+) = Y_{u_k}(\bar{v}^+) = Y_{u_k}(e_1\wedge...\wedge e_k) = -e_1\wedge...\wedge e_{k-1}\wedge e_{2n+1}.\]
Our first goal is to prove the following proposition:

\begin{prop}\label{main prop}
We have $K_k\subseteq V_1(\F)$ and $V_1(\F) \cong V_\F(\lambda_{k-1}^\circ)$, with $v_1^+\in V_1(\F)$ corresponding to a highest weight vector of $V_\F(\lambda_{k-1}^\circ)$.
\end{prop}

\subsection{Proof of Proposition \ref{main prop}}

In the sequel we shall often go back and forth between $V_\F(\lambda_k^\circ)$ and $V_{\bC}(\lambda_k^\circ)$. By a little abuse, we will use the same symbols for vectors of $V_\F(\lambda_k^\circ)$ and $V_{\bC}(\lambda_k^\circ)$, thus using $v^+$ to denote both the highest weight vector chosen in $V_\F(\lambda^\circ_k)$ and the one of $V_{\bC}(\lambda_k^\circ)$. Accordingly, the symbol $v_1^+$, which we have introduced to denote the vector $Y_{u_k}(v^+)$ of $V_\F(\lambda^\circ_k)$ also denotes the vector $Y_{u_k}(v^+)$ of $V_{\bC}(\lambda^\circ_k)$. We do the same in $W_{k}({\bC})$. Thus, the symbol $\bar{v}_1^+ = e_1\wedge...\wedge e_{k-1}\wedge e_{2n+1}$ denotes $\varphi_{k,\bC}(v_1^+)$ as well as $\varphi_{k,\bC}(v_1^+)$, according to the context.  We need one more convention. We put
\[\begin{array}{llll}
2\cdot\mathfrak{A}_{\mathbb{Z}}^\circ &:= & \{2g ~|~ g\in  \mathfrak{A}^\circ_{\mathbb{Z}}\}, & \\
2\cdot V_{\mathbb{Z}}(\lambda_k^\circ) &:= & \{2v ~|~ v\in  V_{\mathbb{Z}}(\lambda_k^\circ)\} & = 2\cdot\mathfrak{A}^\circ_{\mathbb{Z}}(v^+).
\end{array}\]
We say that two vectors $v, w\in V_{\mathbb{Z}}(\lambda^\circ_k)$ are {\em congruent modulo} $2$ and we write $v \equiv_2 w$ if $v-w\in 2\cdot V_{\mathbb{Z}}(\lambda_k^\circ)$. Clearly, the canonical mapping from $V_{\mathbb{Z}}(\lambda_k^\circ)$ to $V_\F(\lambda_k^\circ)$ maps $2\cdot V_{\mathbb{Z}}(\lambda_k^\circ)$ onto $0$. Thus, if $v\equiv_2 w$ in $V_{\mathbb{Z}}(\lambda_k^\circ)$ then $v$ and $w$ are equal when regarded as vectors of $V_\F(\lambda^\circ_k)$.

\begin{lemma}\label{lemma1}
The vector $v_1^+$ belongs to the weight space $V_\mu$ of $V_\F(\lambda^\circ_k)$ with weight $\mu = \lambda^\circ_{k-1} = \lambda_k^\circ - \sum_{i=k}^n\alpha_k$.
\end{lemma}
\pr
We have $u_k = \sum_{i=k}^n\alpha_i$. Hence $v_1^+\in V_\mu$ with $\mu = \lambda^\circ_k-\sum_{i=k}^n\alpha_i$. According to (\ref{lambda kappa}), $\lambda_k^\circ-\sum_{i=k}^n\alpha_i = \sum_{i=1}^{k-2}i\alpha_i + (k-1)\sum_{i=k-1}^n\alpha_i = \lambda_{k-1}^\circ$. \eop

\begin{lemma}\label{lemma2}
We have $X_\alpha(v_1^+) = 0$ for every positive root $\alpha\in\Phi^+$.
\end{lemma}
\pr
We switch from $V_\F(\lambda_k^\circ)$ to $V_{\bC}(\lambda_k^\circ)$. Since $V_{\bC}(\lambda_k^\circ)\cong W_{k}^\circ({\bC})$, we can make our computations inside $W_{k}^\circ({\bC})$. It is easy to see that $X_{u_i-u_j}(\bar{v}_1^+) = X_{u_i+u_j}(\bar{v}_1^+) = 0$ for every choice of $1\leq i < j \leq n$ while $X_{u_i}({\bar v}_1^+) = -2e_1\wedge...\wedge e_{k-1}\wedge e_{i}$ for every $i = 1, 2,..., n$. Therefore $X_{u_i}(v_1^+)\equiv_2 0$ for every
$i = 1, 2,..., n$. It follows that in $V_\F(\lambda_k^\circ)$ we have $X_\alpha(v_1^+) = 0$ for every $\alpha\in \Phi^+$. \eop

\begin{co}\label{corollary3}
The $\mathfrak{A}^\circ_\F$-module $V_1(\F)$ is a homomorphic image of $V_\F(\lambda_{k-1}^\circ)$ and we can choose the homomorphism $\psi_{k-1,\F}:V_\F(\lambda_{k-1}^\circ)\rightarrow V_1(\F)$ in such a way that $\psi_{k-1,\F}$ maps the highest weight vector chosen in $V_\F(\lambda^\circ_{k-1})$ onto $v_1^+$.
\end{co}
\pr
This follows from Theorem \ref{immagini} and Lemmas \ref{lemma1} and \ref{lemma2}. \eop

\bigskip

Our next step is to prove that $K_k\subseteq V_1(\F)$. We firstly compute $Y_{u_i}(v^+)$ and $Y_{u_i}Y_{u_j}(v^+)$ for $i < j$. Recall that $Y_{u_k}(v^+) = v_1^+$ by definition.

\begin{lemma}\label{lemma4}
If $i > k$ then $Y_{u_i}(v^+) = 0$. If $i <  k$ then $Y_{u_i}(v^+) = Y_{u_i-u_k}(v_1^+)$.
\end{lemma}
\pr
We switch to $V_{\mathbb{Z}}(\lambda^\circ_k)$. As $W_{k}(\bC) = W_k^\circ(\bC)\cong V_{\bC}(\lambda^\circ_k)$, we can replace $v^+$ with $\bar{v}^+ = e_1\wedge...\wedge e_k$. Computing in $W_{k}(\bC)$ we see that if $i > k$ then $Y_{u_i}(\bar{v}^+) = 0$ while if $i < k$ then $Y_{u_i}(\bar{v}^+) = Y_{u_i-u_k}Y_{u_k}(\bar{v}^+) = Y_{u_i-u_k}(\bar{v}_1^+)$. Hence in $V_{\mathbb{Z}}(\lambda^\circ_k)$ we have $Y_{u_i}(v^+) = 0$ if $i > k$ and $Y_{u_i}(v^+) = Y_{u_i-u_k}(v_1^+)$ if $i < k$. Clearly, the same holds in $V_\F(\lambda^\circ_k)$. \eop

\begin{co}\label{corollary5}
We have $Y_{u_i}(v^+)\in V_1(\F)$ for every $i = 1, 2,..., n$.
\end{co}

\begin{lemma}\label{lemma6}
If $i < j$ then $Y_{u_i}Y_{u_j}(v^+)\in K_k$.
\end{lemma}
\pr
If $j > k$ then  $Y_{u_i}Y_{u_j}(v^+) = 0$ by Lemma \ref{lemma4}. Let $j = k$. Then $Y_{u_i}Y_{u_j}(\bar{v}^+) =  Y_{u_i}Y_{u_k}(\bar{v}^+) = Y_{u_i}(\bar{v}_1^+) = 2\cdot e_1\wedge...\wedge e_{k-1}\wedge e_{n+i}$, which is $0$ in $W_{k}(\F)$. Therefore $Y_{u_i}Y_{u_k}(v^+)\in K_k$. Finally, let $j < k$.
Then $Y_{u_i}Y_{u_j}(v^+) = Y_{u_i}Y_{u_j-u_k}Y_{u_k}(v^+)$ by Lemma \ref{lemma4}. However $Y_{u_i}$ and $Y_{u_j-u_k}$ commute, since $-u_i-u_j+u_k$ is not a root. Therefore $Y_{u_i}Y_{u_j}(\bar{v}^+) = Y_{u_j-u_k}Y_{u_i}Y_{u_k}(\bar{v}^+) = 0$ since $Y_{u_i}Y_{u_k}(\bar{v}^+) = 0$, as shown above. \eop

\bigskip

We need to fix a few more conventions. The first convention has to do with the order in which we take the factors $Y_\alpha^r/r!$ of a monomial element of type $(-)$ of $\mathfrak{A}^\circ_\F$ (or $\mathfrak{A}^\circ_\mathbb{Z}$). Recall that we only need these elements to produce $V_\F(\lambda^\circ_k)$ (respectively $V_{\mathbb{Z}}(\lambda^\circ_k)$) starting from $v^+$. We assume to write each such element $Y$ as a product $Y = \overline{Y}Y'$ where $Y' = Y_{u_{i_1}}Y_{u_{i_2}}...Y_{u_{i_r}}$ with $i_1 < i_2 <... < i_r$ (and possibly $r = 0$, namely $Y' = 1$) and $\overline{Y}$ is a product of elements as $\frac{Y_{u_i}^2}{2}$, $Y_{u_i+u_j}$, $\frac{Y_{u_i+u_j}^2}{2}$, $Y_{u_i-u_j}$ and $\frac{Y_{u_i-u_j}^2}{2}$ (and possibly $\overline{Y} = 1$), taken in some given order. Note that there is no need to consider elements as $\frac{Y_\alpha^r}{r!}$ with $r > 2$, since $Y_\alpha$ acts on $V_{\mathbb{Z}}(\lambda^\circ_k)$ as a nilpotent endomorphisms of exponent $3$. (Indeed, as noticed in Subsection \ref{orthogonal algebra}, this is true for the action of $Y_\alpha$ on $W_{k}(\bC)$ and $W_{k}(\bC)\cong V_{\bC}(\lambda_k^\circ)$.)

With $Y = \overline{Y}Y'$ as above, we call $\overline{Y}$ (respectively $Y'$) the {\em pseudo-symplectic part} ({\em non-symplectic part}) of $Y$. If $Y' = Y_{u_{i_1}}Y_{u_{i_2}}...Y_{u_{i_r}}$ we put $\delta(Y) := r$ and we call $\delta(Y)$ the {\em non-symplectic degree} of $Y$. If $\delta(Y) = 0$, namely $Y' = 1$, then we say that $Y$ is {\em pseudo-symplectic}, otherwise we say that it is {\em non-symplectic}.  The next corollary immediately follows from Corollary \ref{corollary5} and Lemma \ref{lemma6}.

\begin{co}\label{corollary7}
With $Y$ as above, if $\delta(Y) > 0$ then $Y(v^+)\in V_1(\F)$. If $\delta(Y) > 1$ then $Y(v^+)\in K_k$.
\end{co}

When $\delta(Y) = 1$ a sharper statement can be proved. As in Subsection \ref{symplectic algebra}, let $\overline{V}$ be the hyperplane of $V$ spanned by $\overline{E} = \{e_1,..., e_{2n}\}$ and put $\overline{W}_{k-1} := \bigwedge^{k-1}\overline{V}$. It follows from the information gathered in Subsections \ref{orthogonal algebra} and \ref{symplectic algebra} that $\overline{V}$ is stabilized by all monomials of pseudo-symplectic type. Hence these monomials stabilize $\overline{W}_{k-1}$ too.

By Lemma \ref{lemma4} we immediately obtain the following:

\begin{co}\label{corollary8}
Let $Y$ be a monomial element of $(-)$-type with $\delta(Y) = 1$. Thus $Y = \overline{Y}Y_{u_i}$, where $\overline{Y}$ is the pseudo-symplectic part of $Y$. Define $\bar{w}\in\overline{W}_{k-1}$ as follows: $\bar{w} = 0$ if $i > k$, $\bar{w} = e_1\wedge...\wedge e_{k-1}$ if $i = k$ and
\[\bar{w} = Y_{u_i-u_k}(e_1\wedge...\wedge e_{k-1}) = e_1\wedge... \wedge e_{i-1}\wedge e_k\wedge e_{i+1}\wedge...\wedge e_{k-1}\]
if $i < k$. Then $Y(\bar{v}^+) = \overline{Y}(\bar{w})\wedge e_{2n+1} \in \overline{W}_{k-1}\wedge e_{2n+1}$.
\end{co}

Suppose now that for some monomials $Y_1,..., Y_m$ of $(-)$-type and scalars $t_1,..., t_m\in \F$ we have
\begin{equation}\label{eq0}
(\sum_{i=1}^mt_iY_i)(\bar{v}^+) = 0.
\end{equation}
Regarded $\F$ as a vector space over its prime subfield $\F_2$, let $\tau_0, \tau_1, \tau_2,...$ form a basis of $\F$ over $\F_2$. Thus, $t_i = \tau_0t_{0,i} + \tau_1t_{1,i} + \tau_2t_{2,i}+... $ for suitable scalars $t_{0,i}, t_{1,i}, t_{2,i},...\in \F_2$. Accordingly, (\ref{eq0}) is equivalent to a (possibly infinite) system of equations
\[(\sum_{i=1}^mt_{j,i}Y_i)(\bar{v}^+) = 0, ~~ (j = 0, 1, 2,...).\]
Thus there is no loss in assuming that in (\ref{eq0}) the coefficients $t_1, t_2,..., t_m$ are taken from $\F_2$.

\begin{lemma}\label{lemma9}
Let $g = \sum_{i=1}^mt_iY_i$ be a linear combination of monomials of type $(-)$ with $t_1, t_2,..., t_m\in \F_2$. Suppose that $g(\bar{v}^+) = 0$. Then $g(v^+)\in V_1(\F)$.
\end{lemma}
\pr
We can split $g$ as a sum $g = g_0 + g_1 + g_2$ where $g_0$ is a sum of pseudo-symplectic monomials, $g_1$ is a sum of non-symplectic monomials of degree 1 and $g_2$ is a sum of non-symplectic monomials of degree at least 2. We have $(g_1+g_2)(v^+)\in V_1(\F)$ by Lemma \ref{lemma4}. It remains to prove that $g_0(v^+)\in V_1(\F)$. We shall prove more than this.  Indeed we shall show that $g_0(v^+) = 0$.

We have $g_2(\bar{v}^+) = 0$ by Lemma \ref{lemma6}. Moreover, $g_0(\bar{v}^+)\in \overline{W}_{k}$ (since all monomial elements of pseudo-symplectic type stabilize $\overline{V}$). Let $g_1 = t_rY_r + t_{r+1}Y_{r+1}+... + t_sY_s$ where $Y_j = \overline{Y}_jY_{u_{i_j}}$ for $j = r, r+1,..., s$, $\overline{Y}_j$ being the pseudo-symplectic part of $Y_j$. Put $\widehat{Y}_j = \overline{Y}_jT_j$ where $T_j = 0$ if $i_j > k$, $T_j = 1$ if $i_j = k$ and $T_j = Y_{u_{i_j}-u_k}$ if $i_j < k$. Let $\hat{g}_1 = \sum_{j=r}^st_j\widehat{Y}_j$. By Corollary \ref{corollary8}, we have $g_1(\bar{v}^+) =  (\hat{g}_1(e_1\wedge...\wedge e_{k-1}))\wedge e_{2n+1}$. Therefore
\[g(\bar{v}^+) = g_0(\bar{v}^+) + (\hat{g}_1(e_1\wedge...\wedge e_{k-1}))\wedge e_{2n+1} = 0.\]
It follows that $g(\bar{v}^+) = 0$ if and only if $g_0(\bar{v}^+) = 0$ and $\hat{g}_1(e_1\wedge...\wedge e_{k-1}) = 0$.

Thus $g_0(\bar{v}^+) = 0$ (as $g(\bar{v}^+) = 0$ by assumption) and, without loss of generality, we can assume that $g_1 = g_2 = 0$.  So, $g_0 = g = t_1Y_1 + t_2Y_2 + ... + t_mY_m$ with  $\delta(Y) = 0$ for $i = 1, 2,..., m$. The scalars $t_1, t_2,..., t_m$ belong to $\F_2$. So, $g$ can also be regarded as an element of $\mathfrak{A}^\circ_{\mathbb{Z}}$.

In the sequel of this proof we will freely go back and forth between $\F$ and $\mathbb{C}$ or $\F$ and $\mathbb{Z}$. In order to avoid any confusion, we slightly modify our notation, writing  $v^+_\F$ and $\bar{v}^+_\F$ instead of $v^+$ and $\bar{v}^+$ if we work in $V_\F(\lambda^\circ_k)$ and $W^\circ_{k}(\F)$ and using the symbols $v^+_{\bC}$ and $\bar{v}^+_{\bC}$ when we deal with $V_{\bC}(\lambda^\circ_k)$ or $V_{\mathbb{Z}}(\lambda^\circ)$ and $W^\circ_{k}(\bC)$. However we will do so only for the duration of this proof. Afterwards we will turn back to our earlier simpler notation.

Let us denote by $\overline{\mathfrak{A}}^{\mathrm{sp}}_{\mathbb{Z}}$ the subring of $\mathfrak{A}^{\mathrm{sp}}_{\mathbb{Z}}$ generated by the monomials of $(-)$-type. As we have seen in Subsection \ref{symplectic algebra}, for every monomial $Y$ of pseudo-symplectic type there is a unique monomial $\sigma(Y)\in \overline{\mathfrak{A}}_{\mathbb{Z}}^{\mathrm{sp}}$ such that
\begin{equation}\label{Y-sigma(Y)}
Y(\bar{v}^+_\bC) - \sigma(Y)(\bar{v}^+_\bC) \in 2\cdot\mathfrak{A}_{\mathbb{Z}}^\circ({\bf v}^+_\bC).
\end{equation}
 Explicitly,
\[\sigma(Y_{u_i-u_j}) = V_{u_i-u_j}, \hspace{5 mm}  \sigma(\frac{Y_{u_i}^2}{2}) = V_{2u_i}, \hspace{5 mm}
\sigma(Y_{u_i+u_j}) = V_{u_i+u_j}.\]
We extend $\sigma$ to a product $Y_1Y_2...Y_r$ of monomials $Y_1, Y_2,..., Y_r$ as above by setting $\sigma(Y_1Y_2...Y_r) = \sigma(Y_1)\sigma(Y_2)...\sigma(Y_r)$ and finally we extend it to an  integral combinations $t_1Y_1+ t_2Y_2+...+t_mY_m$ of pseudo-symplectic monomials $Y_1, Y_2,..., Y_m$ by setting $\sigma(\sum_{i=1}^mt_iY_i) = \sum_{i=1}^mt_i\sigma(Y_i)$. Thus, $\sigma(g) = \sum_{i=1}^mt_i\sigma(Y_i)$ and
\begin{equation}\label{g-sigma(g)}
g(\bar{v}^+_\bC)-\sigma(g)(\bar{v}^+_\bC)\in 2\cdot\mathfrak{A}_{\mathbb{Z}}^\circ({\bf v}^+_\bC).
\end{equation}
Turning back to $\F$, since by assumption $g(\bar{v}^+_\F) = 0$ in $W^\circ_{k}(\F)$ we also have $\sigma(g)(\bar{v}^+_\F) = 0$ in $\overline{W}^{\mathrm{sp}}_{k}(\F)$. However, $\overline{W}^{\mathrm{sp}}_{k}(\F) \cong V_\F(\lambda^{\mathrm{sp}}_k)$. Hence $\sigma(g)(u^+_{\F}) = 0$ in $V_\F(\lambda^{\mathrm{sp}}_k)$, where $u_{\F}^+$ stands for the heighest weight vector of $V_\F(\lambda^{\mathrm{sp}}_k)$ corresponding to $\bar{v}^+_\F$ in the isomorphism from $V_\F(\lambda^{\mathrm{sp}}_k)$ to $\overline{W}^{\mathrm{sp}}_{k}(\F)$. This implies that $\sigma(g)(u^+_{\bC})\in 2\cdot V_{\mathbb{Z}}(\lambda^{\mathrm{sp}}_k)$. Whence $\sigma(g)(\bar{v}^+_\bC)\in 2\cdot\overline{\mathfrak{A}}^{\mathrm{sp}}_{\mathbb{Z}}(\bar{v}^+_\bC)$. However, for every element $f\in
\overline{\mathfrak{A}}^{\mathrm{sp}}_{\mathbb{Z}}$ there is an integral combination $h$ of pseudo-symplectic monomial elements of type $(-)$ of $\mathfrak{A}^\circ_{\mathbb{Z}}$ such that $\sigma(h) = f$. We have $h(\bar{v}^+)-f(\bar{v}^+)\in 2\cdot\mathfrak{A}_{\mathbb{Z}}^\circ({\bf v}^+_\bC)$ by (\ref{Y-sigma(Y)}). Therefore $2\cdot\overline{\mathfrak{A}}^{\mathrm{sp}}_{\mathbb{Z}}(\bar{v}^+_\bC)\subseteq 2\cdot\mathfrak{A}^\circ_{\mathbb{Z}}(\bar{v}^+_\bC)$. It follows that
\begin{equation}\label{sigma(g)=0}
\sigma(g)(\bar{v}^+_\bC) \in 2\cdot\mathfrak{A}_{\mathbb{Z}}^\circ({\bf v}^+_\bC).
\end{equation}
By comparing (\ref{g-sigma(g)}) with (\ref{sigma(g)=0}) we obtain that
\begin{equation}\label{g-in-2-I}
g(\bar{v}^+_\bC) \in 2\cdot\mathfrak{A}_{\mathbb{Z}}^\circ({\bf v}^+_\bC).
\end{equation}
Claim (\ref{g-in-2-I}) holds in $W_{k}^\circ(\bC)$. However $W_{k}^\circ(\bC)\cong V_{\bC}(\lambda_k^\circ)$.
Therefore
\begin{equation}\label{g-in-2-II}
g(v^+_\bC) \in 2\cdot\mathfrak{A}_{\mathbb{Z}}^\circ(v^+_\bC) = 2\cdot V_{\mathbb{Z}}(\lambda^\circ_k).
\end{equation}
Namely $g(v^+_\bC) \equiv_2 0$. Claim (\ref{g-in-2-II}) implies that $g(v^+_\F) = 0$, as we wished to prove. \eop

\begin{co}\label{corollary10}
We have $K_k\subseteq V_1(\F)$.
\end{co}
\pr
This immediately follows from Lemma \ref{lemma9}. \eop

\bigskip

Regarded $\mathfrak{A}^\circ_\F$ as an $\F$-vector space, let $\overline{\mathfrak{A}}^\circ_\F$ be its subspace spanned by the set of pseudo-symplectic monomials of $(-)$-type- Then $\overline{\mathfrak{A}}^\circ_\F(\bar{v}^+)$ is a subspace (but not a submodule) of $W^\circ_{k}(\F)$. The following is implicit in the proof of Lemma \ref{lemma9}.

\begin{prop}\label{prop11}
We have $W^\circ_{k}(\F) = \overline{\mathfrak{A}}^\circ_{\F}(\bar{v}^+)\oplus \overline{\mathfrak{A}}^\circ_\F(\bar{v}_1^+)$ (direct sum of vector spaces rather than submodules, since $\overline{\mathfrak{A}}^\circ_{\F}(\bar{v}^+)$ is a subspace but not a submodule of $W^\circ_{k}(\F)$). Moreover $\overline{\mathfrak{A}}^\circ(\bar{v}^+) = \mathfrak{A}^{\mathrm{sp}}_\F(\bar{v}^+) = \overline{W}^{\mathrm{sp}}_{k}(\F)$ (which, as noticed above, is not a submodule of $W^\circ_{k}(\F)$) and $\overline{\mathfrak{A}}^\circ_\F(\bar{v}_1^+) = \mathfrak{A}^\circ_\F(\bar{v}_1^+)$ (which in fact is a submodule of $W^\circ_{k}(\F)$).
\end{prop}

\begin{co}\label{corollary12}
We have $\mathrm{dim}(V_1(\F)) = {{2n+1}\choose{k-1}}$.
\end{co}
\pr
Since $\mathrm{dim}(\overline{W}^{\mathrm{sp}}_{k}(\F)) = {{2n}\choose k}-{{2n}\choose{k-2}}$, Proposition \ref{prop11} implies that $\mathfrak{A}^\circ_\F\bar{v}_1^+$ has codimension ${{2n}\choose k}-{{2n}\choose{k-2}}$ in $W^\circ_{k}(\F)$. On the other hand, $K_k\subset V_1(\F)$ by Corollary \ref{corollary10} and $V_1(\F)/K_k\cong\mathfrak{A}^\circ_\F(\bar{v}_1^+)$. Hence $V_1(\F)$ has codimension ${{2n}\choose k}-{{2n}\choose{k-2}}$. We have ${{2n}\choose k}-{{2n}\choose{k-2}} =
 {{2n+1}\choose k}-{{2n+1}\choose{k-1}}$ in $V_\F(\lambda_k^\circ)$. Hence $\mathrm{dim}(V_1(\F)) = {{2n+1}\choose{k-1}}$, since $\mathrm{dim}(V_\F(\lambda^\circ_k)) = {{2n+1}\choose k}$. \eop

\bigskip

The next corollary finishes the proof of Proposition \ref{main prop}

\begin{co}\label{corollary13}
The homomorphism  $\psi_{k-1,\F}:V_\F(\lambda_{k-1}^\circ)\rightarrow V_1(\F)$ of Corollary {\rm \ref{corollary3}} is indeed an isomorphism.
\end{co}
\pr
Clear from Corollary \ref{corollary12}, since $\mathrm{dim}(V_\F(\lambda
^\circ_{k-1})) = {{2n+1}\choose{k-1}}$.\eop

\bigskip

\noindent
{\bf Remark.} We warn that $\overline{\mathfrak{A}}^\circ_\F$ is not a subalgebra of $\mathfrak{A}^\circ_\F$. For instance, $[Y_{u_i-u_j},\frac{Y_{u_j}^2}{2}] = Y_{u_i+u_j}-Y_{u_j}Y_{u_i} = -Y_{u_i+u_j}-Y_{u_i}Y_{u_j}$. Each of $Y_{u_i-u_j}$, $\frac{Y_{u_j}^2}{2}$ and $Y_{u_i+u_j}$ belongs to $\overline{\mathfrak{A}}^\circ_\F$ but the sum $Y_{u_i+u_j}-Y_{u_j}Y_{u_i}$ does not. Therefore $[Y_{u_i-u_j},\frac{Y_{u_j}^2}{2}]\not\in \overline{\mathfrak{A}}^\circ_\F$. In other words, only one of the products $Y_{u_i-u_j}\frac{Y_{u_j}^2}{2}$ and $\frac{Y_{u_j}^2}{2}Y_{u_i-u_j}$ belongs to $\overline{\mathfrak{A}}^\circ_\F$, depending on the order chosen for the factors occurring in the pseudo-symplectic monomial elements.

On the other hand, $Y_{u_i}Y_{u_j}(\bar{v}^+) = 0$ by Corollary \ref{corollary7}. Hence $[Y_{u_i-u_j},\frac{Y_{u_j}^2}{2}](\bar{v}^+) = -Y_{u_i+u_j}(\bar{v}^+)$. Thus, if for instance $\overline{\mathfrak{A}}^\circ_\F$ contains the product $Y_{u_i-u_j}\frac{Y_{u_j}^2}{2}$ (whence $\frac{Y_{u_j}^2}{2}Y_{u_i-u_j}\not\in \overline{\mathfrak{A}}^\circ_\F$), then the space $\overline{\mathfrak{A}}^\circ_\F(\bar{v}^+)$ also contains the vector $\frac{Y_{u_i}^2}{2}Y_{u_i-u_j}(\bar{v}^+)$, but represented as a sum
$\frac{Y_{u_i}^2}{2}Y_{u_i-u_j}(\bar{v}^+) = Y_{u_i-u_j}\frac{Y_{u_j}^2}{2}(\bar{v}^+)+Y_{u_i+u_j}(\bar{v}^+)$.

\subsection{End of the proof of Theorem \ref{main theorem}}

In order to finish the proof of Theorem \ref{main theorem} it remains to prove that $V_1(\F) = {\cal N}_k$ (notation as in Section \ref{Introduction}) and that the isomorphism $\psi_{k-1,\F}:V_\F(\lambda_{k-1}^\circ)\rightarrow V_1(\F)$ maps the subspace ${\cal N}_{k-1}$ of $V_\F(\lambda_{k-1}^\circ)$ onto $K_k$.

\begin{prop}\label{prop14}
We have $V_1(\F) = {\cal N}_k$.
\end{prop}
\pr
As $K_k\subseteq V_1(\F)$, in order to prove that $V_1(\F) = {\cal N}_k$ we only must prove that $V_1(\F)/K_k = {\cal N}_k/K_k$, namely the homomorphism $\varphi_{k,\F}:V_\F(\lambda_k^\circ)\rightarrow W^\circ_{k}(\F)$ maps $V_1(\F)$ onto the $k$-nucleus subspace $\overline{\cal N}_k$ of $W^\circ_{k}(\F)$.

We recall that $\overline{\cal N}_k = \langle \iota_{k}(N_X)\rangle_{X\in \Delta_{k-1}}$, where $\Delta_{k-1}$ is the set of totally singular $(k-1)$-subspaces of $V$, $N_0$ is the nucleus of the quadric described by $\eta$ and $N_X = \langle X, N_0\rangle$ for $X\in \Delta_{k-1}$.

According to the conventions stated at the beginning of Subsection \ref{orthogonal algebra}, we have $N_0 = \langle e_{2n+1}\rangle$. Put $X_1 := \langle e_1, e_2,..., e_{k-1}\rangle$. Then $X_1\in \Delta_{k-1}$ and $N_{X_1} = \langle e_1,..., e_{k-1}, e_{2n+1}\rangle$. Accordingly, $\iota_{k}(N_{X_1}) = \langle e_1\wedge...\wedge e_{k-1}\wedge e_{2n+1}\rangle = \langle \bar{v}_1^+\rangle$. Therefore $\varphi_{k,\F}(V_1(\F))$ contains $\iota_{k}(N_{X_1})$. As $\varphi_{k,\F}(\F)$ is a $G$-module and $G$ acts transitively on $\Delta_{k-1}$, the module $\varphi_{k,\F}(V_1(\F))$ contains $\overline{\cal N}_k$.

On the other hand, both $\overline{\cal N}_k$ and $\varphi_{k,\F}(V_1(\F))$ have codimension ${{2n+1}\choose k}-{{2n+1}\choose{k-1}}$ in $W^\circ_{k}(\F)$ (see Theorem \ref{theorem 1} and Corollary \ref{corollary12}). Hence $\varphi_{k,\F}(V_1(\F)) = \overline{\cal N}_k$.  \eop

\bigskip

\noindent
{\bf Remark.} The argument used in the proof of Proposition \ref{prop14} can be exploited to prove
that $\mathrm{dim}(V_1(\F)) = {{2n+1}\choose{k-1}}$ (Corollary \ref{corollary12}) avoiding Proposition \ref{prop11}. Indeed by that argument we obtain that $\varphi_{k,\F}(V_1(\F))\supseteq \overline{\cal N}_k$. Therefore $\mathrm{dim}(V_1(\F)) \geq {{2n+1}\choose{k-1}}$, since $\overline{\cal N}_k$ has codimension ${{2n+1}\choose k}-{{2n+1}\choose{k-1}}$ in $W^\circ_{k}(\F)$. On the other hand,
 $\mathrm{dim}(V_1(\F)) \leq {{2n+1}\choose{k-1}}$ as $V_1(\F)$ is a homomorphic image of $V_\F(\lambda^\circ_{k-1})$. Hence $\mathrm{dim}(V_1(\F)) = {{2n+1}\choose{k-1}}$.

\begin{prop}\label{prop15}
The isomorphism $\psi_{k-1,\F}:V_\F(\lambda_{k-1}^\circ)\rightarrow V_1(\F)$ maps the subspace ${\cal N}_{k-1}$ of $V_\F(\lambda_{k-1}^\circ)$ onto $K_k$.
\end{prop}
\pr
Put $v_2^+ := Y_{k-1}v_1^+$ and $V_2(\F) := \mathfrak{A}^\circ_\F(v_2^+)$. We must prove that $V_2(\F) = K_k$.  We have $Y_{k-1}(\bar{v}_1^+) = Y_{k-1}Y_k(\bar{v}^+) = 0$ (by Corollary \ref{corollary7}). Therefore $V_2(\F)\subseteq K_k$.

The pre-images $\psi_{k-1,\F}^{-1}(V_2(\F))$ and $\psi_{k-1,\F}^{-1}(v_1^+)$ of $V_2(\F)$ and $v_1^+$ are the analogues of $V_1(\F)$ and $v_0^+$ respectively, but in $V_\F(\lambda_{k-1}^\circ)$ instead of $V_\F(\lambda^\circ_k)$. So, we can apply Proposition \ref{main prop} to them, obtaining that $\psi_{k-1,\F}^{-1}(V_2(\F))\cong V_\F(\lambda^\circ_{k-2})$, whence
$V_2(\F)\cong V_\F(\lambda^\circ_{k-2})$. It follows that $\mathrm{dim}(V_2(\F)) = {{2n+1}\choose{k-2}}$. However we also have $\mathrm{dim}(K_k) = {{2n+1}\choose{k-2}}$ (see Section \ref{Introduction}, (\ref{dim-grass})). Hence $V_2(\F) = K_k$, since $V_2(\F)\subseteq K_k$. \eop

\bigskip

The proof of Theorem \ref{main theorem} is complete.

\section{Uniqueness conjectures}\label{uniqueness}

So far we have proved that $V_\F(\lambda_k^\circ)$ admits a series of submodules
\[0 \subset M_0 \subset M_1 \subset .... \subset M_{k-1}\subset M_k = V_\F(\lambda^\circ_k)\]
where
\begin{equation}\label{chain0}
M_i = Y_{u_{i+1}}Y_{u_{i+2}}....Y_{u_{k-1}}Y_{u_k}(v^+) \cong V_\F(\lambda_i^\circ) , ~~ (i = 0, 1, 2,..., k-1).
\end{equation}
Moreover $M_i$ corresponds to the submodule ${\cal N}_{i+1}$ of $V_\F(\lambda_{i+1}^\circ)$ as well as to the kernel $K_{i+2}$ of the homomorphism $\varphi_{i+2,\F}:V_\F(\lambda_{i+2}^\circ)\rightarrow W_{i+2}^\circ(\F)$ (if $i < k-1$). Avoiding any mention of the elements
$Y_{u_{i+1}}....Y_{u_k}(v^+)$ we can rephrase the above conditions as follows:
\begin{equation}\label{chain1}
M_i \cong V_\F(\lambda_i^\circ) ~~ (i = 0, 1, 2,..., k),
\end{equation}
\begin{equation}\label{chain2}
M_i/M_{i-2} \cong W^\circ_{i}(\F) ~~ (i = 1, 2,..., k).
\end{equation}
Recall that $V_\F(\lambda_0^\circ)$, to be considered in (\ref{chain1}) for $i = 0$, is $1$-dimensional, by convention. We allow $i = 1$ in (\ref{chain2}). This forces us to define $M_{-1}$ too. We put $M_{-1} := 0$.

When $\F$ is perfect we also know that
\begin{equation}\label{chain3}
M_{i}/M_{i-1} \cong \overline{W}^{\mathrm{sp}}_{i}(\F).
\end{equation}
(Theorem \ref{theorem 1 bis} and Corollary \ref{main corollary ter}.)  The following conjecture is quite natural.

\begin{conj}\label{conj1}
The series of submodules defined as in {\rm (\ref{chain0})} is the unique series of $G$-submodules of $V_\F(\lambda^\circ_k)$ satisfying conditions {\rm (\ref{chain1})} and {\rm (\ref{chain2})}.
\end{conj}

If the $G$-module $V_\F(\lambda^\circ_k)$ was rigid, namely all of its automorphisms as a $G$-module fix $\langle v^+\rangle$, then Conjecture \ref{conj1} would be not so difficult to prove, but we do not know if $V_\F(\lambda^\circ_k)$ is rigid (although it is certainly rigid as an $\mathfrak{A}^\circ_\F$-module). Recall that all irreducible Weyl modules are rigid (Humphreys \cite[2.2]{Humphreys2}), but when $\mathrm{char}(\F) = 2$, as we are assuming here, $V_\F(\lambda^\circ_k)$ is not irreducible.

When $\F$ is perfect we can consider the following variation of Conjecture \ref{conj1}.

\begin{conj}\label{conj2}
Let $\F$ be perfect. Then the series of submodules defined as in {\rm (\ref{chain0})} is the unique series of $G$-submodules of $V_\F(\lambda^\circ_k)$ satisfying conditions {\rm (\ref{chain1})} and {\rm (\ref{chain3})}.
\end{conj}

We shall prove Conjecture \ref{conj2} for $2 \leq k \leq 4$ (when $k = 1$ both Conjectures \ref{conj1} and \ref{conj2} are trivial). We will rely on the information on the lattice of submodules of $\overline{W}^{\mathrm{sp}}_{k}(\F)$ provided in Baranov and Suprunenko \cite{BS} (see also Premet and Suprunenko \cite{PS}). As the reader will see, in principle our method can be applied to check Conjecture \ref{conj2} for every particular value of $k$, but apparently there is no way to make a general argument out of it.

Let $\F$ be a perfect field of characteristic $2$. Let $(M_i)_{i=0}^k$ be the series defined as in (\ref{chain0}) and let $(M'_i)_{i=0}^k$ be another series  satisfying conditions (\ref{chain1}) and (\ref{chain3}).

\bigskip

\noindent
1) Let $k = 2$ and suppose by contradiction that $M'_1 \neq M_1$. Then either $M_1\cap M'_1 = M_0$ or $M_1\cap M'_1 = 0$. Hence $(M_1+M'_1)/M_1$ occurs in $M_2/M_1\cong \overline{W}^{\mathrm{sp}}_{2}(\F)$ as a submodule isomorphic to either $\overline{W}^{\mathrm{sp}}_{1}(\F)$ or $W^\circ_{1}(\F)$. But according to \cite{BS} and \cite{PS} no such submodules occur in $\overline{W}^{\mathrm{sp}}_{2}(\F)$.

Therefore $M'_1 = M_1$. As $M'_1 = M_1 \cong V_\F(\lambda^\circ_1)$ admits a unique $1$-dimenional submodule, we also have $M'_0 = M_0$.

\bigskip

\noindent
2)  Let $k = 3$. According to \cite{BS} and \cite{PS}, if $n$ is odd then the module $M_2/M_1 \cong \overline{W}^{\mathrm{sp}}_{2}(\F)$ is irreducible while if $n$ is even then it admits a unique proper submodule $\overline{M}_{2,1}$, of dimension $1$. In this case let $M_{2,1}\supset M_1$ be such that $M_{2,1}/M_1 = \overline{M}_{2,1}$. As for $M_1 \cong W_{1}^\circ(\F)$, we know that $M_0$ is its unique proper submodule. The module $W^\circ_{2}(\F)$ does not admit any $1$-dimensional submodule, as one can check by direct computations. Hence $M_{2,1}/M_0$ does not split as a direct sum of $M_2/M_0$ and a $1$-dimensional module.

It follows from the above that if $M'_2\neq M_2$ the following are the only possibilities for $M_2\cap M'_2$:
\[\begin{array}{ccll}
M_2\cap M'_2 & \mathrm{dim}(M_2\cap M'_2) & \mathrm{dim}((M_2+M'_2)/M_2) & \\
{} & & & \\
0 & 0 &  {{2n+1}\choose 2} &  \\
M_0 & 1 &  {{2n+1}\choose 2}-1 & \\
M_1 & 2n+1 & {{2n+1}\choose 2}-2n-1 & \\
M_{2,1} & 2n+2 & {{2n+1}\choose 2}-2n-2 &  (\mbox{only if $n$ is even})
\end{array}\]
On the other hand, according to \cite{BS} and \cite{PS}, if $n$ is even then the module $\overline{W}^{\mathrm{sp}}_{3}(\F)$ is irreducible while if $n$ is odd then it admits a unique proper submodule, isomorphic to $\overline{W}^{\mathrm{sp}}_{1}(\F)$. Therefore ${{2n}\choose 3}-2n$ and $2n$ are the only possibilities for $\mathrm{dim}((M_2+M'_2)/M_2)$. Neither of the numbers listed in the second column of the above table is equal to any of these two numbers, except for ${{2n+1}\choose 2}-2n-1$, which is equal to ${{2n}\choose 3}-2n$ for $n = 3$.

So, we are left with the case of $n = 3$, $M'_2\cap M_2 = M_1$ and $M'_2+M_2 = M_3 = V_\F(\lambda_3^\circ)$. However $W^\circ_{3}(\F) = V_\F(\lambda_3^\circ)/K_3 = M_3/M_1$. Whence $W^\circ_{3}(\F)$ now splits as the direct sum of two copies $M_2/M_1$ and $M'_2/M_1$ of $\overline{W}^{\mathrm{sp}}_{2,3}(\F)$. (Note that, in view of the analysis made in the case $k = 2$, $M_1 = M'_1$ is the unique submodule of $M'_2$ of dimension $2n+1 = 7$.) Accordingly, the stabilizer in $G$ of the $1$-dimensional subspace spanned by $\bar{v}_1^+ = Y_{u_3}\bar{v}^+\in M_2/M_1$ in $G$ also stabilizes a $1$-dimensional subspace of $M'_2/M_1$. It is not difficult to see by direct computations that this is not the case.

It follows that $M'_2 = M_2$. The equalities $M'_1 = M_1$ and $M'_0 = M_0$ now follow by the analysis of case $k = 2$ applied to $M_2 = V_\F(\lambda_2^\circ)$.

\bigskip

\noindent
3) Let $k = 4$. As seen in the previous cases, if $n$ is odd then $M_0$ and $M_1$ are the unique non-zero proper submodules of $M_2$ while if $n$ is even then $M_2$ also admits a submodule $M_{2,1}\supset M_1$ with $\mathrm{dim}(M_{2,1}/M_1) = 1$. We have also seen that $M_3/M_2$ is irreducible if $n$ is even while if $n$ is odd then $M_3/M_2$ admits a unique proper submodule $\overline{M}_{3,1}$, isomorphic to $\overline{W}^{\mathrm{sp}}_{1}(\F)$. Let $M_{3,1}\supset M_2$ be such that $M_{3,1}/M_2 = \overline{M}_{3,1}$. By considering the action of $G$ on $M_3/M_1 = W^\circ_{3}(\F)$ it is not difficult to see that $M_{3,1}/M_1$ cannot split as a direct sum $M_{3,1}/M_1 = M_2/M_1\oplus \overline{M}_{3,1}$ of $G$-submodules. Therefore the following are the unique non-zero proper $G$-submodules of $M_3$:
\[\begin{array}{ccc}
\mbox{$n$ even} & \hspace{5 mm} & \mbox{$n$ odd}\\
{} & & \\
\begin{array}{cc}
\mbox{submodule} & \mbox{dimension}
{}{}\\
M_0 & 1 \\
M_1 & 2n+1 \\
M_{2,1} & 2n+2 \\
M_2 & {{2n+1}\choose 2} \\
M_3 & {{2n+1}\choose 3}
\end{array} & &
\begin{array}{cc}
\mbox{submodule} & \mbox{dimension}
{}{}\\
M_0 & 1 \\
M_1 & 2n+1 \\
M_2 & {{2n+1}\choose 2} \\
M_{3,1} & {{2n+1}\choose 2}+2n\\
M_3 & {{2n+1}\choose 3}
\end{array}
\end{array}\]
Accordingly, the following are the only possibilities for $M_3\cap M'_3$ if $M'_3\neq M_3$:
\[\begin{array}{ccll}
M_3\cap M'_3  & \hspace{5 mm} & \mathrm{dim}((M_3+M'_3)/M_3) &  \\
{} & & &  \\
0 & & {{2n+1}\choose 3} &  \\
M_0 & & {{2n+1}\choose 3}-1 & \\
M_1 & & {{2n+1}\choose 3}-2n-1 & \\
M_{2,1} & & {{2n+1}\choose 3}-2n-2 &  (\mbox{only if $n$ is even}) \\
M_2 & & {{2n+1}\choose 3}-{{2n+1}\choose 2} & \\
M_{3,1} & & {{2n+1}\choose 3}-{{2n+1}\choose 2}-2n & (\mbox{only if $n$ is odd})
\end{array}\]
According to \cite{BS} and \cite{PS}, if $n \equiv 3 ~(\mathrm{mod}~4)$ then the module $\overline{W}^{\mathrm{sp}}_{4}(\F)$ is irreducible while if $n \equiv 1~ (\mathrm{mod}~4)$ then it admits a unique non-zero proper submodule, which is $1$-dimensional. If $n \equiv 0 ~(\mathrm{mod}~4)$ then $\overline{W}^{\mathrm{sp}}_{4}(\F)$ admits a unique non-zero proper submodule, which has dimension equal to ${{2n}\choose 2}-2$ and is isomorphic to the unique irreducible quotient of $\overline{W}^{\mathrm{sp}}_{2}(\F)$. Finally, if $n \equiv 2~ (\mathrm{mod}~4)$ then $\overline{W}^{\mathrm{sp}}_{4}(\F)$ admits just two proper submodules, namely either $\overline{W}^{\mathrm{sp}}_{2}(\F)$ and its unique $1$-dimensional submodule or the dual of $\overline{W}^{\mathrm{sp}}_{2}(\F)$ and its unique submodule of codimension $1$. Therefore
the following are the only possibilities for $\mathrm{dim}((M_3+M'_3)/M_3)$.
\[\begin{array}{cl}
1 & \mbox{(only if $n\equiv 1 ~(\mathrm{mod}~4)$ or $n\equiv 2~ (\mathrm{mod}~4)$)}\\
{{2n}\choose 2}-2 & \mbox{(only if $n$ is even)}\\
{{2n}\choose 2}-1 &  \mbox{(only if $n\equiv 2~ (\mathrm{mod}~4)$)}\\
{{2n}\choose 4}-{{2n}\choose 2} &
\end{array}\]
Comparing these values with the feasible dimensions previously listed for $(M_3+M'_3)/M_3$ we see that
$M_3\cap M'_3 = 0$ with $n = 5$ is the unique possiblity. In this case $(M_3+M'_3)/M_3 = M_4/M_3$. This implies that $M_4/M_2\cong W^\circ_{4}(\F)$ splits as the direct sum of two $G$-submodules, one of which is $M_3/M_2\cong \overline{W}^{\mathrm{sp}}_{3}(\F)$ and the other one is $(M'_3+M_2)/M_2\cong M'_3\cong V_\F(\lambda_3^\circ)$. By an argument similar to the one used in a similar situation in case $k = 3$, one can see that $W^\circ_{4}(\F)$ does not admit such a splitting.

Therefore $M'_3 = M_3$. We can now apply to $M_3\cong V_\F(\lambda_3^\circ)$ the analysis made in case $k = 3$, obtaining that $M'_2 = M_2$, $M'_1 = M_1$ and $M'_0 = M_0$.

\bigskip

\noindent
{\bf Remark.}  An analysis as above has also been done in \cite{CP2}, but only for $k \leq n\leq 3$.

\bigskip

\noindent
Authors' address\\

\noindent
Ilaria Cardinali, Antonio Pasini,\\
Department of Information Engineering and Mathematics,\\
University of Siena,\\
Via Roma 56, 53100 Siena, Italy\\
ilaria.cardinali@unisi.it, antonio.pasini@unisi.it

\end{document}